\documentclass[journal]{IEEEtran}
\usepackage{cite}
\usepackage{url}
\usepackage{empheq}
\usepackage{float}
\usepackage{dsfont}
\usepackage[hidelinks]{hyperref}
\usepackage[normalem]{ulem}
\usepackage{amsmath,amssymb,amsfonts}
\allowdisplaybreaks[4]
\usepackage{graphicx}
\usepackage{textcomp}
\usepackage{amsmath}
\usepackage{enumerate}
\usepackage{epstopdf}
\usepackage{array}
\usepackage{booktabs}
\usepackage{subfigure}
\usepackage{multirow}
\usepackage{soul, color}
\usepackage[usenames,dvipsnames]{xcolor}
\usepackage[latin1]{inputenc}
\usepackage{cases}

\newtheorem{theorem}{Theorem}
\newtheorem{lemma}{Lemma}

\newtheorem{proposition}{Proposition}

\soulregister\cite7 
\soulregister\citep7 
\soulregister\citet7 
\soulregister\ref7 
\soulregister\pageref7

\usepackage{bbding}

\usepackage{nomencl}
\makenomenclature
\usepackage{etoolbox}
\renewcommand\nomgroup[1]{%
  \item[\bfseries
  \ifstrequal{#1}{A}{Acronyms}{%
  \ifstrequal{#1}{S}{Symbols}{%
  \ifstrequal{#1}{U}{Units}{}}}%
]}

\usepackage{amsmath}
\allowdisplaybreaks[4]

\usepackage{algorithm}
\usepackage[noend]{algpseudocode}

\begin{document}

\title{Carbon-Aware Quantification of Real-Time Aggregate Power Flexibility of Electric Vehicles}

\author{Xiaowei~Wang,~\IEEEmembership{Student Member,~IEEE}, and
Yue~Chen,~\IEEEmembership{Senior Member,~IEEE}
\thanks{This work was supported by the National Natural Science Foundation of China under Grant No. 52307144 and the 1+1+1 CUHK-CUHK(SZ)-GDST Joint Collaboration Fund General R\&D Project. (\emph{Corresponding to Y. Chen})}
\thanks{X. Wang and Y. Chen are with the Department of Mechanical and Automation Engineering, the Chinese University of Hong Kong, HKSAR, China. E-mail:  \{xwwang, yuechen\}@mae.cuhk.edu.hk.}
}
\maketitle

\begin{abstract}
Electric vehicles (EVs) can be aggregated to offer flexibility to power systems. However, the rapid growth in EV adoption leads to increased grid-level carbon emissions due to higher EV charging demand, challenging grid decarbonization efforts. Quantifying and managing the EV flexibility while controlling carbon emissions is crucial. This paper introduces a methodology for carbon-aware quantification of real-time aggregate EV power flexibility based on the Lyapunov optimization technique. 
We construct a novel queue system including EV charging queues, delay-aware virtual queues, and carbon-aware virtual queues. Based on the evolution of these queues, we define the Lyapunov drift and minimize the drift-plus-penalty term to get the real-time EV flexibility interval, which is reported to the system operator for flexibility provision.
To enhance EV flexibility, we integrate dispatch signals from the system operator into the queue updates through a two-stage disaggregation process. The proposed approach is prediction-free and adaptable to various uncertainties, including EV arrivals and grid carbon intensity. Additionally, the maximum charging delay of EV charging tasks is theoretically bounded by a constant, and carbon emissions are effectively controlled. The numerical results demonstrate the effectiveness of the proposed online method and highlight its advantages over several benchmarks through comparisons.
\end{abstract}


\IEEEpeerreviewmaketitle
\section*{Nomenclature}
\addcontentsline{toc}{section}{Nomenclature}
\begin{IEEEdescription}[\IEEEusemathlabelsep\IEEEsetlabelwidth{$Q(p,r^{\pm},w)$}]
\item[$t,\mathcal{T}$] Index and set of time periods.
\item[$i,\mathcal{I}$] Index and set of EVs.
\item[$k,\mathcal{K}$] Index and set of EV groups.
\item[$t_i^a,t_i^d$] Arrival time and departure time of EV $i$.
\item[$e_i^{ini},e_i^{req}$] Initial and required energy level of EV $i$.
\item[$R_i$] The allowed charging duration of EV $i$.
\item[$\check{p}_{s,t},\hat{p}_{s,t}$] Lower and upper bounds of aggregate EV power flexibility interval at time $t$.
\item[$\Delta t$] Time interval.
\item[$\delta_c$] Charging efficiency.
\item[$\check{p}_{i,t},\hat{p}_{i,t}$] Lower and upper bounds of EV $i$'s charging power at time $t$.
\item [$p_{i,max}$] Charging power limit of EV $i$.
\item[$\check{e}_{i,t},\hat{e}_{i,t}$] Lower and upper bounds of EV $i$ 's battery energy at time $t$.
\item[${e}_{i,min},{e}_{i,max}$] Minimum and maximum energy of EV $i$.
\item [$w_{g,t}$] Carbon intensity of grid electricity at time $t$.
\item [$E_s$] Carbon emission quota of the charging station over the operation horizon.
\item [$J_{t}$] EV charging queue backlog at time $t$.
\item [$a_t$] Arrival rate of EV charging tasks at time $t$.
\item [$a_{i,t}$] Arrival rate of the charging task generated by EV $i$ at time $t$.
\item [$J_{k,t}$] EV charging queue backlog of group $k$ at time $t$.
\item[$\check{p}_{k,t},\hat{p}_{k,t}$] Lower and upper bounds of EV power flexibility interval of group $k$ at time $t$.
\item [$a_{k,t}$] Arrival rate of EV charging tasks of group $k$ at time $t$.
\item [$H_{k,t}$] Delay queue backlog of group $k$ at time $t$.
\item [$\lambda$] Parameter associated with delay queue.
\item [$R_k$] Charging duration of group $k$.
\item [$r$] Time-average carbon emission rate cap over the operation horizon.
\item [$Q_{c,t}$] Carbon queue backlog at time $t$.
\item [$\Bar{P}_{k,t}$] Charging power limit of group $k$ at time $t$.
\item [$\Bar{P}_{i,t}$] Maximum charging power of EV $i$ at time $t$.
\item [$\beta$] Parameter associated with carbon queue.
\item [$V$] Weight parameter in the objective function.
\item [$\gamma_t$] Dispatch ratio at time $t$.
\item [$p_{s,t}^d$] Dispatch signal at time $t$.
\item [$p_{k,t}^d$] Dispatch signal for group $k$ at time $t$.
\item [$p_{1,k,t}^d$] Dispatch signal for group $k$ in the first stage at time $t$.
\item [$p_{2,k,t}^d$] Dispatch signal for group $k$ in the second stage at time $t$.
\end{IEEEdescription}

\section{Introduction}
\IEEEPARstart{T}{he} power system is experiencing a growing penetration of fluctuating renewable energy sources due to the energy crisis and climate challenges, which greatly challenges the system security. Utilizing the demand side flexibility is a promising way to address these challenges \cite{8635327}. In particular, the growing electric vehicles (EVs) are one of such demand-side flexibility assets, whose sizable energy storage and controllable load-shifting capability enable EV charging scheduling to be coordinated into the system energy management without compromising the comfort of EV owners \cite{10180215}. While the decarbonization of the transportation sector has accelerated thanks to the growing adoption of EVs, the rising charging demand of massive EVs could lead to an increase in the carbon emissions in electricity generation \cite{DIXON20201072}. Therefore, it is crucial to quantify the EV flexibility with consideration of carbon emission control to facilitate an effective and environmentally friendly power system.

For easier integration into the system-level energy management, it is desired to quantify the aggregate flexibility of massive EVs, which could effectively react to the dispatch signal from the system operator\cite{10286155}. Some works have been devoted to the quantification of aggregate EV power flexibility. The summation of power and energy bounds of all EVs was adopted to evaluate the aggregate EV flexibility, which was then utilized to improve the economic benefits \cite{EVagg_sum1} and the reliability\cite{EVagg_sum2} of power systems. Reference \cite{EVagg_battery} employed a virtual battery model to depict the aggregate EV flexibility in alleviating network congestion. The aggregation of EV flexibility was formulated by a dispatchable region to explore its usage in microgrid bidding \cite{TIIwang}. An inner approximation method was proposed to maximize the EV flexibility \cite{TSEyan}. Reference \cite{EVagg_cluster} constructed an EV clustering and aggregation method considering the spatio-temporal correlations of EVs. A polytope-based aggregation method was designed for EV aggregation considering EV heterogeneity \cite{TSG_polytope}.

The works above solve the aggregate EV power flexibility characterization problem offline, which requires complete information of future EV arrivals. However, the highly stochastic behaviors of EVs are hard to predict precisely. Moreover, the dispatch signal from the system operator is not considered in the offline model, potentially resulting in an underestimation of the EV flexibility. Therefore, an online (also called real-time) aggregate EV power flexibility quantification method that does not require complete future information and considers the dispatch signal would be more practical. The real-time aggregate EV flexibility is usually defined as an aggregate charging power adjustment range \cite{EVagg_def} and many works have studied the problem of real-time flexibility quantification. The online aggregate flexibility in unbalanced systems was characterized by the inner approximation method in a model predictive control (MPC) manner \cite{chen2019aggregate}. Based on the MPC, the real-time flexibility of residential buildings was described by a flexibility envelope in \cite{GASSER2021116653} and a flexibility band in \cite{9147459}. However, the MPC framework relies on predictions to perform the rolling process, and the forecasting inaccuracy may decrease the performance. Different from the MPC methods, data-driven approaches are promising alternatives which are conducted based on real-time observations. Reinforcement learning (RL) was used to solve the real-time aggregate EV flexibility with feedback design \cite{RL_li}.
Reference \cite{RL_zhang} proposed an RL-based approach to estimate the real-time aggregate EV flexibility. To address the issue of high computational burdens in real-time flexibility quantification, AdaBoost regressor was proposed for online use \cite{10202703}. Nevertheless, the learning approaches rely on a large amount of data for training and cannot ensure adaptability to different situations due to the limitation of the quality of datasets.

Although the aforementioned methods demonstrate promising results in real-time aggregate power flexibility quantification, they do not take carbon emission control into account in their designs.
Reference \cite{9960988} formulated a carbon-aware EV charging problem with the objective of minimizing carbon emissions. A carbon-efficient real-time feedback system for EV charging was proposed in \cite{HUBER2021124766}. 
However, the above two works relied on grid carbon density forecasts and focused on optimizing the real-time EV charging rather than quantifying the EV flexibility.
The carbon-based online energy flexibility of heat pumps was tested in \cite{PEAN2019101579}, and the carbon-aware real-time flexibility of buildings was explored in \cite{CAI2024105531}. However, both references were based on the MPC, requiring perfect predictions of the future. 

Lyapunov optimization is a promising method in real-time applications, where decisions are made under uncertain environments without future predictions \cite{neely2010stochastic}. Lyapunov optimization has been utilized in various applications in smart grid, including microgrid control {\cite{Lya0}, \cite{Lya1}}, battery storage control \cite{Lya4}, and real-time energy management of data centers \cite{Lya5}. Particularly, in references \cite{Lya2, Lya3}, queue systems were constructed to store the EV charging requests, and Lyapunov optimization was used to minimize the charging cost of EVs under uncertainty. Moreover, an online EV charging scheduling algorithm based on Lyapunov optimization was proposed to alleviate the network congestion in~\cite{mohan2024market}. Although Lyapunov optimization performs well in real-time applications, applying it to the real-time aggregate EV power flexibility quantification problem with consideration of carbon emission control is more challenging in three aspects: 1) Cost efficiency does not necessarily mean carbon efficiency. Reducing emissions based on time-varying carbon intensity may delay EV charging, while prioritizing charging can compromise emission targets. 2) The current dispatch signal from the power system operator may impact the future EV flexibility, which is hard to take into account. 3) Maximizing EV flexibility while ensuring their carbon footprints within the emission limit is difficult.

We compare the related works in TABLE \ref{tab:literature} to clarify the research gaps. This paper aims to fill the research gaps by proposing a carbon-aware online algorithm based on Lyapunov optimization for real-time aggregate EV power flexibility quantification to address the research challenges. The main contributions are summarized as follows: 

1) \emph{Modeling}:
For online applications, we model the aggregate EV power flexibility using time-decoupled flexibility intervals. To formulate the carbon-aware EV flexibility quantification problem for solving these flexibility intervals, we construct a novel queue system, including EV charging queues for handling real-time charging requests, delay-aware virtual queues for bounding the worst-case EV charging delays, and carbon-aware virtual queues for controlling the carbon emission rate of the system. With these queues, the flexibility quantification problem is formulated as a stochastic programming problem, aiming to maximize the time-average total aggregate EV power flexibility subject to time-average constraints on EV charging and carbon emission, to fit within the framework of Lyapunov optimization.

2) \emph{Algorithm}: 
The Lyapunov optimization technique is employed to calculate the EV flexibility intervals in real time by minimizing the drift-plus-penalty term. Compared to the traditional Lyapunov optimization methods, the proposed method solves the problem by minimizing the quadratic form of the drift-plus-penalty term rather than its linear approximation. This results in a more uniformly distributed aggregate flexibility region with an enhanced total flexibility value. To pass the real-time dispatch signal from the system operator to individual EVs, we propose a two-stage disaggregation process. In the first stage, we address backlogs in the EV charging queue to complete charging tasks as quickly as possible. In the second stage, the unallocated power after the first stage is distributed according to the EV arrival times. The dispatch signal is also used to update all queue backlogs, accurately reflecting the states of the queues, and thereby enhancing EV flexibility. Note that the proposed method is prediction-free and theoretically guaranteed to achieve near-optimality in maximizing EV flexibility.

\begin{table*}[t]
\small
    \centering
    \caption{Comparison of Related Online Approaches In Existing Literature} 
    \vspace{-0.5em}
    \begin{tabular}{@{}ccccccccc@{}}
        \toprule
       {Ref.} & {Online method} & {Carbon awareness} & {Flexibility quantification} & {Large-scale EVs} & {No future information}\\
        \midrule
        \cite{chen2019aggregate,GASSER2021116653,9147459} & MPC & & \checkmark & & \\
        \cite{9960988,HUBER2021124766} & MPC & \checkmark &  &  &  \\
         \cite{PEAN2019101579,CAI2024105531} & MPC & \checkmark &\checkmark  &  &  \\
        \cite{RL_li} & Learning method & & \checkmark & \checkmark & \checkmark \\
         \cite{RL_zhang} & Learning method & & \checkmark & \checkmark &  \\
         \cite{10202703} & Learning method & & \checkmark & & \checkmark \\
         \cite{Lya2, Lya3,mohan2024market} & Lyapunov optimization & & &\checkmark & \checkmark \\
        This work  & Lyapunov optimization & \checkmark & \checkmark & \checkmark & \checkmark\\
        
        \bottomrule
        \label{tab:literature}
   \end{tabular}
\end{table*}

The remainder of this paper is organized as follows. Section \ref{sec:formulation} introduces the system framework and formulates the offline model to describe the carbon-aware aggregate EV power flexibility quantification problem. Section \ref{sec:online} establishes the queue system to reformulate the problem into a time-average form and develops an online algorithm based on Lyapunov optimization to solve it for real-time aggregate EV flexibility quantification. Section \ref{sec:Disaggregation} performs a two-stage disaggregation to achieve the dispatch signal, which is used in the subsequent queue update procedure. Numerical simulations are conducted in Section \ref{sec:case}. Section \ref{sec:conclusion} concludes the paper.

\section{Problem Formulation}
\label{sec:formulation}
This section begins by introducing the overall framework for aggregate EV power flexibility quantification and then the offline model considering the carbon emission constraint.

\begin{figure}[!htpb]
    \centering
    \includegraphics[width = 0.85\linewidth]{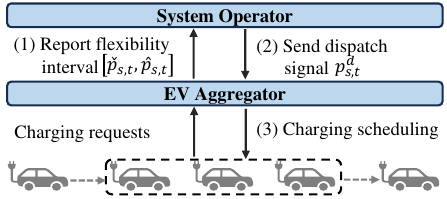} 
    \caption{System diagram for EV power flexibility quantification.}
    \label{fig:schematic}
\end{figure}

\subsection{System Framework}
As shown in Fig. \ref{fig:schematic}, an aggregator determines the charging scheduling of all EVs in a charging station. The EVs communicate with the aggregator in real time to submit their charging requests. $\mathcal{I}$ is the set of EVs and the time slot is indexed by $t \in \mathcal{T} := \{1,..., T\} $. When a vehicle $i \in \mathcal{I}$ arrives at the charging station, the EV owner will inform the aggregator of his or her charging task, described by four parameters ($t_i^a, t_i^d, e_i^{ini},e_i^{req}$), where $t_i^a$ is the arrival time of EV $i$, $t_i^d$ is the departure time, $e_i^{ini}$ is its initial battery energy level upon arrival and $e_i^{req}$ is the desired energy level when leaving. The allowed/maximum charging duration of EV $i$ is $R_i = t_i^d - t_i^a$, during which the EV $i$ can be dynamically adjusted to different charging rates as long as the charging task is completed before departure, which provides a certain degree of EV power flexibility. By collecting the information of all in-station EVs, the aggregator can quantify the real-time aggregate EV power flexibility of the charging station.

Generally, more power consumption leads to more carbon emissions. Therefore, more power flexibility provision would also result in more carbon emissions. It is challenging to quantify the real-time aggregate EV power flexibility with carbon emission constraints since the EV charging behaviors are time-coupled and their carbon emission footprints have time-accumulated effects. The above issues motivate us to depict the aggregate EV power flexibility by time-decoupled regions for online use, such that any aggregate charging power within this region can be disaggregated into feasible charging power of individual EVs without violating both EV charging and carbon emission constraints. Therefore, in each time slot, the aggregator quantifies a flexibility interval $[\Check{p}_{s,t},\hat{p}_{s,t}]$ and reports it to the system operator to support grid-level energy management. This flexibility intervals are described by their lower and upper power bounds (i.e., $\Check{p}_{s,t}$ and $\hat{p}_{s,t}$) in each time slot and can be easily integrated into the system operation with little communication load. Upon receiving the flexibility interval, the system operator can choose any power level within the flexibility interval and send a dispatch signal $p_{s,t}^d \in [\Check{p}_{s,t},\hat{p}_{s,t}]$ to the aggregator to utilize this EV flexibility. The aggregator then disaggregates the dispatch signal into the charging power for individual EVs.

\subsection{Offline Carbon-Aware EV Power Flexibility Quantification}
In this section, we first give an offline formulation to illustrate the carbon-aware aggregate EV power flexibility quantification problem. Then we list several drawbacks regarding this offline problem to better motivate our transition to an online setting.

The offline algorithm is assumed to have perfect information about future EV arrivals and the time-varying carbon intensity of the electricity from the grid.  A set of time-decoupled flexibility intervals $[\Check{p}_{s,1},\hat{p}_{s,1}] \times ... \times [\Check{p}_{s,T},\hat{p}_{s,T}]$ are used to describe the aggregate EV flexibility region. The upper and lower bounds of the flexibility region are specified by $\{\hat{p}_{s,t}, \forall{t} \}$ and $\{\check{p}_{s,t} \forall{t} \}$, respectively. In the following, we continue to use $\check{(\cdot)}$ and $\hat{(\cdot)}$ to denote the variables associated with the lower and upper power trajectories.  Let $\Delta t$ be the duration of one time slot. Then we derive this region by maximizing the energy level of the aggregate EV flexibility over the optimization horizon, i.e.,
\begin{equation}\label{eq:off_obj}
\mathbf{P1\colon} \quad \max \quad  \sum\nolimits_{t \in \mathcal{T}} \left[(\Hat{p}_{s,t} -\Check{p}_{s,t})\Delta t - \epsilon (\Hat{p}_{s,t} -\Check{p}_{s,t})^2\right].
\end{equation}
Note that the problem $\mathbf{P1}$ may have multiple solutions. To get a unique solution, a quadratic term with a small enough constant $\epsilon > 0 $ (i.e., $- \epsilon (\Hat{p}_{s,t} -\Check{p}_{s,t})^2$) is added to the objective function \eqref{eq:off_obj} to allocate the flexibility intervals across the time horizon as evenly as possible.

The constraints of the upper trajectory are given by:
\begin{subequations}\label{eq:off_ub}
    \begin{align}
        & \Hat{p}_{s,t} = \sum\nolimits_{i \in \mathcal{I}} \hat{p}_{i,t}, \forall{t}, \label{eq:off_ub_1}\\
&0 \leq \hat{p}_{i,t} \leq p_{i,max}, \forall{i},\forall{t} \in [t_i^a,t_i^d],\label{eq:off_ub_2}\\
& \hat{p}_{i,t}= 0,\forall{i},\forall{t} \notin [t_i^a,t_i^d],\label{eq:off_ub_3}\\ 
& \hat{e}_{i,t_i^a} = e_i^{ini}, \hat{e}_{i,t_i^d} \geq e_i^{req}, \forall{i}, \label{eq:off_ub_4}\\
& \hat{e}_{i,t} = \hat{e}_{i,t-1} + \delta_c\hat{p}_{i,t-1}\Delta t, \forall{i},\forall t \in \mathcal{T}/\{1\},\label{eq:off_ub_5}\\
& e_{i,min}\leq \hat{e}_{i,t} \leq e_{i,max}, \forall{i}, \forall t,\label{eq:off_ub_6}
    \end{align}
\end{subequations}
where \eqref{eq:off_ub_1} shows that the aggregate power $ \hat{p}_{s,t}$ is the summation of the charging power of individual EVs $\hat{p}_{i,t},\forall i $. \eqref{eq:off_ub_2} shows that the charging power of EV $i$ is constrained by its charging power limit $p_{i,max}$. \eqref{eq:off_ub_3} ensures that the EV can only be charged during the charging session. Let $e_{i,t}$ denote the battery energy level of EV $i$ at time $t$. Constraint \eqref{eq:off_ub_4} shows the initial energy level of the EV and its required energy when leaving the charging station. Constraint \eqref{eq:off_ub_5} gives the state-of-charge (SoC) dynamics and $\delta_c$ is the charging efficiency. \eqref{eq:off_ub_6} sets the upper and lower bounds of EV $i$'s battery energy level, denoted by $e_{i,max}$ and $e_{i,min}$, respectively.

Similarly, we can have the following constraints for the corresponding lower power trajectory:
\begin{subequations}\label{eq:off_lb}
    \begin{align}
        &  \check{p}_{s,t} = \sum\nolimits_{i \in \mathcal{I}} \check{p}_{i,t}, \forall{t}, \label{eq:off_lb_1}\\ 
&0 \leq \check{p}_{i,t} \leq p_{i,max}, \forall{i},\forall{t} \in [t_i^a,t_i^d],\label{eq:off_lb_2}\\ 
& \check{p}_{i,t}= 0,\forall{i},\forall{t} \notin [t_i^a,t_i^d],\label{eq:off_lb_3}\\ 
& \check{e}_{i,t_i^a} = e_i^{ini}, \check{e}_{i,t_i^d} \geq e_i^{req}, \forall{i}, \label{eq:off_lb_4}\\
& \check{e}_{i,t} = \check{e}_{i,t-1} + \delta_c\check{p}_{i,t-1}\Delta t, \forall{i}, \forall t \in \mathcal{T}/\{1\},\label{eq:off_lb_5}\\ 
& e_{i,min}\leq \check{e}_{i,t} \leq e_{i,max}, \forall{i}, \forall t. \label{eq:off_lb_6}
    \end{align}
\end{subequations}

The following joint constraint ensures that the power of the upper trajectory should be no less than the lower trajectory.
\begin{equation} \label{eq:ub_lb}
 \qquad \; \hat{p}_{s,t} \geq \Check{p}_{s,t}, \forall{t}.
\end{equation}

The EV charging demand of massive EVs leads to substantial electricity consumption, resulting in significant carbon emission, as it serves as the driver of electricity generation that relies on carbon-intensive fuels.
Therefore, a carbon emission cap $E_s$ is imposed on the upper trajectory so that the carbon emission of any trajectory of the flexibility region over the optimization horizon is restricted to be no larger than the given carbon emission cap. 
\begin{equation}\label{eq:off_carbon}
    \quad\sum \nolimits_{t \in \mathcal{T}} w_{g,t}\Hat{p}_{s,t} \Delta t \leq E_s.
\end{equation}

Note that for charging stations owned by the power grid, the power grid determines such carbon emission cap to comply with the grid decarbonization target. For other stations, such as commercial charging stations, the government determines the value of $E_s$ to respond to climate change, which depends on the regional carbon emission requirements and policies \cite{10807052}. In \eqref{eq:off_carbon}, $w_{g,t}$ is the carbon intensity of the electricity from the grid at time $t$, which is sensitive to both location and time \cite{10.1145/3632775.3661970}. 

\begin{figure}[!htpb]
    \centering
    \includegraphics[width = 0.95\linewidth]{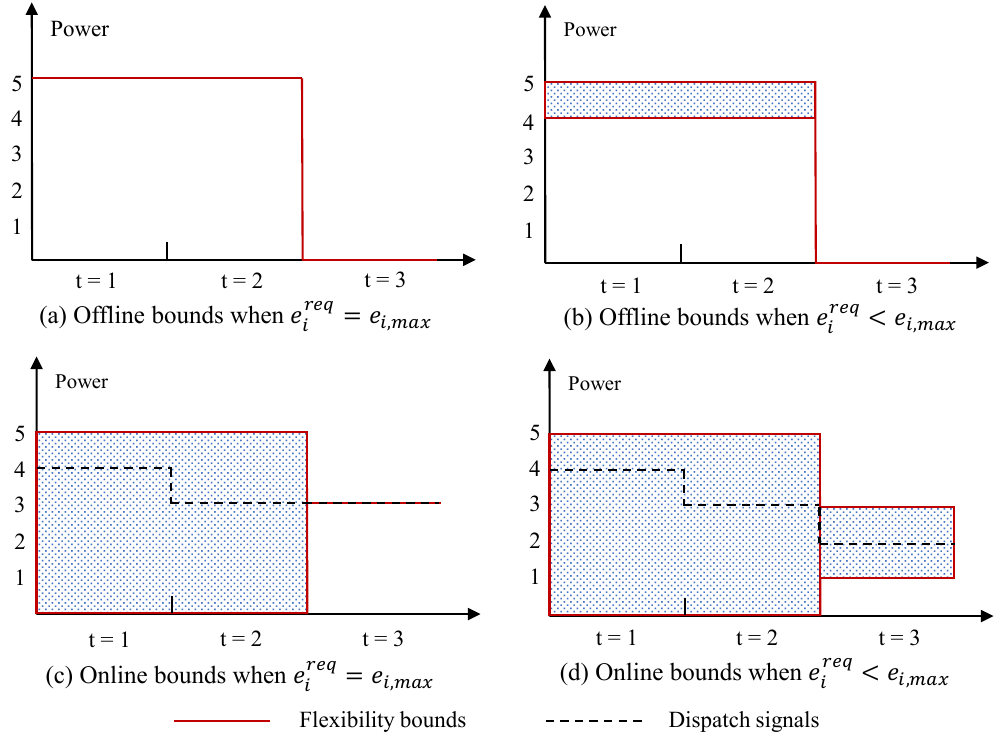} 
    \caption{A three-time slot example showing the 
flexibility regions obtained by the offline and online methods.}
    \label{fig:offline_online_cmp}
\end{figure}

The offline model has the following three drawbacks: 1) It requires perfect information of future EV charging behaviors and time-varying grid carbon intensity over all periods, which is hard to know accurately. 2) The model cannot reflect the impact of the real-time dispatch signal issued by the system operator on the future EV flexibility, leading to an underestimation of the EV flexibility. 3) Due to constraints \eqref{eq:off_ub_4} and \eqref{eq:off_lb_4}, the model may fail to capture EV flexibility. This is because when the required energy of the EV is equal to the upper limit of its battery energy (i.e., $e_{i}^{req} = e_{i,max}$), the upper and lower bounds of the flexibility interval coincide.

To tackle the above limitations, we develop an online method to derive the real-time aggregate EV power flexibility based on Lyapunov optimization. 
In an online setting, we can utilize real-time observations to handle the uncertainties in the problem. For uncertainties related to EV charging requests, we can update the number of EVs available at the station and their corresponding charging demands with the real-time EV arrival observations. For uncertain grid carbon intensity, we can observe the real-time carbon intensity and integrate this information into the real-time flexibility quantification to meet the carbon constraint. 
Furthermore, in an online setting,
the dispatch signal from the system operator acts as a feedback signal, which can be utilized by the EV aggregator to update the EV battery states. Through this feedback design, the limitation of the offline method when  $e_i^{req}$ is equal to $e_{i,max}$ is overcome.

Fig. \ref{fig:offline_online_cmp} provides an intuitive explanation of how the online method overcomes the limitation of the offline model by using a toy example with three time slots. In this illustrative example, we omit units for simplicity. We assume the initial EV battery energy is zero with a maximum charging rate of 5 units per time slot. The upper limit of EV battery energy $e_{i,max}$ is assumed to be equal to the energy accumulated from two time slots of charging at maximum rate (i.e., 10 units). Fig. \ref{fig:offline_online_cmp}(a) shows that when $e_{i}^{req} = e_{i,max}$, the upper and lower bounds of the offline model for all three time slots overlap, failing to show EV flexibility. However, in Fig. \ref{fig:offline_online_cmp}(c), the online method shows non-overlapping bounds. To be specific, at time slot 1, the remaining parking time of the EV is larger than the least time needed to complete the charging requests. Thus, the online algorithm can adaptively maintain the lower bound at a reduced level to provide EV flexibility. The system operator then issues a dispatch signal (i.e., 4 units per time slot) to the aggregator, which subsequently updates the EVs' battery energy level to 4 units. A similar process is repeated for time slot 2 with a dispatch signal of 3 units per time slot. When the final time slot 3 is reached, the upper and lower bounds are all 3 units (i.e., $10-4-3=3$ units) per time slot. One can observe that the upper bound cannot be further increased due to the constraint of $e_{i,max}$, while the lower bound cannot be further reduced because it shall ensure the fulfillment of the remaining EV charging demands in the final time slot. When $e_{i}^{req} = 8$ units (which is less than $e_{i,max}$), the offline model demonstrates a degree of flexibility since the lower bounds can be less than upper bounds, which is shown  Fig. \ref{fig:offline_online_cmp}(b). For the online method, the flexibility in Fig. \ref{fig:offline_online_cmp}(d) is similarly to  Fig. \ref{fig:offline_online_cmp}(c) for the first two time slots. However, for the third time slot, the online approach can exhibit additional flexibility due to the reduced required energy $e_{i}^{req}$. The difference between the offline and online methods stems from the online approach's ability to dynamically adjust EV battery states by incorporating dispatch signals as real-time feedback.


It is worth exploring how to relate the real-time flexibility interval to EV charging requests and carbon emission constraints.
Intuitively, the upper bounds of the flexibility intervals over the operation horizon cannot remain high, as this would give the system operator more opportunities to dispatch higher power levels, potentially leading to higher carbon emission and thus violating the carbon constraint. The lower bounds cannot be too low for extended periods, as this would result in low charging power to EVs, potentially causing a large number of unfulfilled EV charging requests. 
Moreover, since the system operator can choose any power level within the flexibility interval, we need to ensure that the required constraints (such as the carbon emission constraint) are met when we determine the flexibility intervals during the aggregation process. Once the real-time dispatch signal is received, we can update our online model in the subsequent disaggregation process with this dispatch signal to enhance future EV flexibility. Therefore, the aggregation and disaggregation processes are not independent and can influence each other.

In the next section, to achieve the goals discussed above, we will construct a novel queue system to capture the system dynamics and employ the Lyapunov optimization technique to quantify the real-time aggregate EV power flexibility.

\begin{proposition}\label{pr:1}
For any given aggregate power trajectory $\{p_{s,t}^{reg}\}$ that satisfies $\check{p}^{\ast}_{s,t} \leq p_{s,t}^{reg} \leq \hat{p}^{\ast}_{s,t},\forall t $, there exists an achievable carbon-aware charging strategy.
\end{proposition}

The proof of Proposition \ref{pr:1} is provided in Appendix A. 

\section{Online Algorithm}
\label{sec:online}
In this section, to develop the online counterpart of the offline model, we first construct the queue models to describe the system dynamics.
Then, we formulate the problem as stochastic programming to optimize the time-average aggregate EV power flexibility subject to time-average constraints on both EV charging and carbon emission. Finally, based on the queue backlogs, an online algorithm that does not require future information is designed to solve the problem in real time by Lyapunov optimization with theoretical guarantees.

\subsection{EV Charging Queue Model}

Inspired by the data packets in communication systems, we define the EV charging queue model with dynamics as follows:
\begin{equation}\label{eq:queue_J}
    J_{t+1} = \max[J_{t} - \check{p}_{s,t},0]+a_t,
\end{equation}
where $J_t$ is called the queue backlog at time slot $t$, which can represent the EV charging tasks that need to be fulfilled. Recall that in the offline model $\mathbf{P1}$, the lower bound of the flexibility region is closely related to the energy demand of the EVs, the decision variable $\check{p}_{s,t}$ stands for the amount of charging tasks that the queue can process at time $t$. $a_t$ is the stochastic EV arrival process, whose value represents the new charging tasks that arrive during time slot $t$. The charging tasks are served based on the First-In-First-Out (FIFO) principle. It is worth noting that $a_t$ is first available for processing at time slot $t+1$. $a_t$ is the summation of all individual EV charging tasks arriving during time slot $t$, which is shown as follows:
\begin{equation}
    a_t = \sum\nolimits_{i \in \mathcal{I}}a_{i,t},
\end{equation}
where $a_{i,t}$ is the arriving charging task of EV $i$ during time slot $t$. 
For one time slot, an individual EV $i$ can add at most a charging task of $p_{i,max}$ (i.e., the charging power limit of EV $i$) to the EV charging queue backlog, which is similar to a new data package that needs to be processed by a server in a communication system. In general, multiple time slots are required to finish the charging demand of the EV. Therefore, one EV can add charging tasks to more than one time slot. For example, suppose an EV has a charging power limit of 5kW, and its charging energy demand can be fulfilled within two time slots. We further assume that during the first time slot, the EV is charged at its maximum power limit of 5kW, while in the second time slot, a charging power of 3kW is sufficient to meet its remaining demand. Therefore, the EV initially adds a charging task of 5kW upon connection to the charger and subsequently adds a charging task of 3kW in the next time slot. Note that if the 5 kW demand in the first time slot is not fully processed (for instance, only 4 kW is processed), then the remaining 1 kW will be stored in the backlog. In the next time slot, this remaining 1 kW will be prioritized over the new arriving charging task of 3kW for processing, as we follow the FIFO principle. 
Thus, $a_{i,t}$ is determined by both the charging power limit and the required energy of EV $i$, which is generated by:
\begin{equation}\label{eq:ev_task}
 a_{i,t} = 
\begin{cases}
p_{i,max}, &t_i^a \leq t <  \lfloor{\eta_i}\rfloor + t_i^a, \\
e_i^{task}/(\delta_c\Delta t) - \lfloor{\eta_i}\rfloor p_{i,max},&t = \lfloor{\eta_i}\rfloor + t_i^a,\\
0, &\text{otherwise}.
\end{cases}
\end{equation}
In \eqref{eq:ev_task}, $e_i^{task} = (e_{i}^{req}-e_i^{ini})$ is the required energy of EV $i$, $\eta_i = e_i^{task}/(\delta_c p_{i,max} \Delta t)$ is the required timeslots when the EV is charged at the maximum charging rate, and $\lfloor{\eta_i}\rfloor$ is the greatest integer less than or equal to $\eta_i$.

Since EV owners select different charging duration, the acceptable charging delay of their charging tasks could also be different, which motivates us to define multiple EV charging queues for different charging durations. We construct $K$ queues to store the EV charging tasks, and each queue corresponds to a charging duration of $R_k, k \in \mathcal{K} :=\{1,2,...K\}$. The resulting queue dynamics are driven by:
\begin{equation}\label{eq:queue_J_k}
    J_{k,t+1} = \max[J_{k,t} - \check{p}_{k,t},0]+a_{k,t}, \forall{k},
\end{equation}
where $J_{k,t}$, $\check{p}_{k,t}$ and $a_{k,t}$ are the corresponding queue backlog, EV charging task processing rate, and arriving rate at time slot $t$ of queue $k$. In particular, we have the following:
\begin{equation}
    J_t = \sum_{k \in \mathcal{K}}J_{k,t},~ \check{p}_{s,t} = \sum_{k \in \mathcal{K}}\check{p}_{k,t}, ~a_t = \sum_{k \in \mathcal{K}}a_{k,t}.
\end{equation}

\subsection{Delay-Aware Virtual Queue}
Similar to the delay of data packets in communication systems, the delay of the EV charging tasks also needs to be considered.  Therefore, to ensure that the worst-case charging task delay is bounded, we introduce a delay-aware virtual queue $H_{k,t}$ for each group $k$ with its initial value of backlog $H_{k,1} = 0$, whose update is dependent on $J_{k,t}$ and $\check{p}_{k,t}$.  
 Specifically, if $ J_{k,t} > \check{p}_{k,t}$, the resulting update equation is as follows:
\begin{equation}\label{eq:queue_H_k}
H_{k,t+1} = \max[H_{k,t} + \frac{\lambda}{R_k} - \check{p}_{k,t},0], \forall{k}.
\end{equation}
Otherwise, we set $H_{k,t+1}$ to $0$. In the above equation, 
$\lambda$ is a given parameter. We can see that $\frac{\lambda}{R_k}$ serves as a penalty to the backlog of the queue $H_{k,t}$ when the charging tasks in the queue $J_{k,t}$ is not cleared (where the condition $ J_{k,t} > \check{p}_{k,t}$ is active). Furthermore, the constant $\lambda$ controls the growth rate of queue $H_{k,t}$, and intuitively, a longer charging duration $ R_k$ of group $k$ leads to a smaller $\frac{\lambda}{R_k}$ for queue $H_{k,t}$. This is reasonable since more attention should be paid to the queue with shorter charging duration, which renders queue $H_{k,t}$ delay-aware. The size of queue $H_{k,t}$ provides a bound on the delay of the EV charging tasks in queue $J_{k,t}$. To be specific, for each group $k$, if a scheduling algorithm can guarantee that $J_{k,t}$ and $H_{k,t}$ are upper bounded by finite constants $J_{k,max}$ and $H_{k,max}$ , respectively, for all timeslots, then the worst-case delay of charging tasks is also bounded by $\delta_{k,max}$, which is defined by:
\begin{equation}\label{def_delay}
    \delta_{k,max} \stackrel{\triangle}{=} \frac{(J_{k,max}+H_{k,max})R_k}{\lambda}.
\end{equation}

\begin{lemma}\label{lemma:1}
    Assume $J_{k,t}$ and $H_{k,t}$ evolve with dynamics described by \eqref{eq:queue_J_k} and \eqref{eq:queue_H_k}, respectively, and that $J_{k,t} \leq J_{k,max}$ and $H_{k,t} \leq H_{k,max}$ for all time slots $t \in \mathcal{T}$ can be ensured by an algorithm. Suppose EV charging tasks are serviced in a FIFO order. Then the worst-case delay of all EV charging tasks in queue $k$ is bounded by $\delta_{k,max}$.
\end{lemma}

The proof of Lemma \ref{lemma:1} can be found in Appendix B.

\subsection{Carbon-Aware Virtual Queue}
According to \eqref{eq:off_carbon}, the carbon emission is the time accumulation of carbon emission rate (i.e., $w_{g,t}\hat{p}_{s,t}$), which is determined by the time-varying grid carbon intensity $w_{g,t}$ and the aggregate charging power $\hat{p}_{s,t}$. The carbon emission limit is a long-term consideration. Therefore, we can  decompose this long-term limit into individual per-slot problems to design our online method. We divide both sides of \eqref{eq:off_carbon} by $T$:
\begin{equation}\label{carbon_1}
  \frac{\sum_{t \in \mathcal{T}} w_{g,t}\Hat{p}_{s,t} \Delta t}{T} \leq \frac{E_s}{T}.
\end{equation}
To ensure that the total carbon emission during $T$ time slots does not exceed the budget $E_s$, the time-average carbon emission must not exceed ${E_s}/{T}$. We further divide two sides of \eqref{carbon_1} by $\Delta t$:
\begin{equation}\label{carbon_2}
  \frac{\sum_{t \in \mathcal{T}} w_{g,t}\Hat{p}_{s,t}}{T} \leq r,
\end{equation}
where $r = {E_s}/{(T \Delta t)}$. Considering $ w_{g,t}\Hat{p}_{s,t}$ is the carbon emission rate of time slot $t$, equation \eqref{carbon_2} indicates that we can control the time-average carbon emission rate to ensure the accumulated carbon emission over the optimization horizon will not exceed the carbon emission limit $E_s$. Therefore, parameter $r$ serves as a cap imposed on the time-average carbon emission rate. To achieve the goal specified in \eqref{carbon_2}, we additionally construct the carbon-aware virtual queue $Q_{c,t}$ to capture the time-varying carbon emission rate with dynamics:
\begin{equation}\label{eq:queue_Q_c}
    Q_{c,t+1} = \max[Q_{c,t} + w_{g,t}\hat{p}_{s,t} - r,0],
\end{equation}
where $Q_{c,t}$ is the backlog of the carbon-aware virtual queue, which captures the time accumulated violations of the carbon emission rate constraint \eqref{carbon_2}. A larger backlog of $Q_{c,t}$ reflects that many previous time periods have carbon emission rates exceeding $r$ (i.e., $w_{g,t}\hat{p}_{s,t} - r > 0$), indicating inefficient control of the carbon emission footprint of EV charging.

Note that different from the definition of the EV charging queue in \eqref{eq:queue_J} that uses the lower bound of the flexibility interval $\check{p}_{s,t}$ as the aggregate charging power, the carbon-aware virtual queue $Q_{c,t}$ uses the upper bound of the flexibility interval $\hat{p}_{s,t}$ as the aggregate charging power. This is because, intuitively, the upper bound of the flexibility region is more restricted by the carbon emission constraint. Besides, in \eqref{eq:queue_Q_c}, if we multiply $w_g\hat{p}_{s,t}$ and $r$ by $\Delta t$, then $Q_{c,t}$ can capture the accumulated violations of the absolute carbon emission
constraint \eqref{carbon_1}.
However, as in the previous queue model (e.g., \eqref{eq:queue_J_k}), our goal is to more intuitively illustrate how the carbon constraint impacts aggregate EV power and how the aggregate power relates to the queue backlog. Therefore, we have discarded $\Delta t$ and focused on carbon emission rate control rather than absolute carbon emissions. Based on the derivation above, the two can be equivalently converted into each other.

\subsection{Lyapunov Optimization}
In the following, we will restate the problem in a time-average sense and use the queue systems that we have constructed to meet the required constraints. Then, we apply the Lyapunov optimization technique to solve the problem.

We transform the problem $ \mathbf{P1}$ into the following to maximize the time-average energy level of aggregate EV power flexibility:
\begin{subequations}\label{eq:P2}
    \begin{align}
 \mathbf{P2\colon} &\min_{\check{p}_{k,t},\hat{p}_{k,t},\forall{k,t}}\;   \lim_{T \to \infty}\frac{1}{T}\sum_{t=1}^{T} \mathbb{E}[-(\Hat{p}_{s,t} -\Check{p}_{s,t})\Delta t], \label{eq:P2_a}\\
&\textrm{s.t.} \quad \lim_{T \to \infty}\frac{1}{T}\sum_{t=1}^{T} \mathbb{E}[a_{k,t} - \check{p}_{k,t}] \leq 0, \forall{k},\label{eq:P2_b}\\
& \qquad \lim_{T \to \infty}\frac{1}{T}\sum_{t=1}^{T} \mathbb{E}[ w_{g,t}\hat{p}_{s,t} - r] \leq 0, \label{eq:P2_c}\\
& \qquad \qquad 0 \leq \check{p}_{k,t} \leq \Bar{P}_{k,t}, \forall{k,t},\label{eq:P2_d}\\
& \qquad \qquad 0 \leq \hat{p}_{k,t} \leq \Bar{P}_{k,t}, \forall{k,t},\label{eq:P2_e}\\
& \qquad \qquad \quad \check{p}_{k,t} \leq \hat{p}_{k,t}, \forall{k,t}, \label{eq:P2_f} \\
& \qquad   \hat{p}_{s,t} = \sum_{k \in \mathcal{K}}\hat{p}_{k,t},~\check{p}_{s,t} = \sum_{k \in \mathcal{K}}\check{p}_{k,t}.
    \end{align}
\end{subequations}

In the objective, we use the fact that maximizing the average flexibility level is equivalent to minimizing its negative value. Constraint \eqref{eq:P2_b} sets a limit that for group $k$, the lower bound $\check{p}_{k,t}$ shall meet the charging demand $a_{k,t} $ in a time-average sense. \eqref{eq:P2_c} ensures that the time-average carbon emission rate of the upper bound $\hat{p}_{s,t}$ should not exceed the parameter $r$. \eqref{eq:P2_d} and \eqref{eq:P2_e} show that in each group, both the lower and the upper bounds of the aggregate power are restricted by the maximum charging power of group $k$, which is denoted by $\Bar{P}_{k,t}$ and satisfies $\Bar{P}_{k,t} = \sum_{i \in \mathcal{I}_k}\Bar{P}_{i,t}$, where $\Bar{P}_{i,t}, i \in \mathcal{I}_k$ is the maximum charging power of EV $i$ in group $k$ at time $t$ and can be calculated by the following:
\begin{align}
    \Bar{P}_{i,t} = \left\{
    \begin{aligned}
        & P_{i,max}, \begin{aligned}
            ~\\~
        \end{aligned} \\
        & \frac{e_{i,max} - e_{i,t}}{\delta_c \Delta t}, \begin{aligned}
            ~\\~
        \end{aligned} \\
        & 0, 
    \end{aligned}~~~~
    \begin{aligned}
        & \begin{aligned}
            & \text{if} \;t_i^a \leq t <  t_i^d \;\text{and}\\
            & e_{i,t} + \delta_c p_{i,max}\Delta t\leq e_{i,max},
        \end{aligned} \\
        & \begin{aligned}
            & \text{if} \;t_i^a \leq t <  t_i^d \;\text{and}\\
            & e_{i,t} + \delta_c p_{i,max}\Delta t > e_{i,max},
        \end{aligned} \\
        & \text{otherwise}.
    \end{aligned}
    \right.
\end{align}

Constraint \eqref{eq:P2_f} ensures that the upper bound of the flexibility interval in group $k,\forall{k}$ is greater than the lower bound.

Next, we show how queue stability can help meet the time-average constraints \eqref{eq:P2_b} and \eqref{eq:P2_c}. By the definition of $J_{k,t}$ in \eqref{eq:queue_J_k}, for each $k \in \mathcal{K}$, we have:
\begin{equation}\label{eq:queue_J_k_inequ}
    J_{k,t+1} - a_{k,t} \geq J_{k,t} - \check{p}_{k,t}, \forall{t}.
\end{equation}
Summing the above inequality over all time slots $t \in \{1,2,...,T\}$ ,dividing two sides by $T$ and taking the expectations, then we have:
\begin{equation}
    \frac{\mathbb{E}[J_{k,T+1}] - \mathbb{E}[J_{k,1}] }{T} \geq \frac{\sum_{t=1}^{T} \mathbb{E}[a_{k,t} - \check{p}_{k,t}]}{T}. 
\end{equation}
When $\lim_{T \to \infty} \mathbb{E}[J_{k,T+1}]/T = 0 $ holds for the backlog of EV charging queue $J_{k,t}$, which means that if the queue $J_{k,t}$ is mean-rate-stable, then constraint \eqref{eq:P2_b} is satisfied. Following the same procedure, we can make the carbon-aware virtual queue $Q_{c,t}$ mean-rate-stable (i.e., $\lim_{T \to \infty}\mathbb{E}[Q_{c,T+1}]/{T} = 0$) to meet the constraint \eqref{eq:P2_c}.

Replacing constraints \eqref{eq:P2_b} and \eqref{eq:P2_c} with queue stability constraints, we can obtain a problem as follows:
\begin{equation}\label{eq:2}
    \begin{aligned}
 \mathbf{P3\colon} \min_{\check{p}_{k,t},\hat{p}_{k,t},\forall{k,t}}\;   &\lim_{T \to \infty}\frac{1}{T}\sum_{t=1}^{T} \mathbb{E}[-(\Hat{p}_{s,t} -\Check{p}_{s,t})\Delta t],\\
 \qquad &\textrm{s.t.} \quad (\text{16d}) - (\text{16g}),\\
 &\lim_{T \to \infty} \frac{\mathbb{E}[J_{k,T+1}]}{T} = 0, \forall{k},\\
 &\lim_{T \to \infty} \frac{\mathbb{E}[Q_{c,T+1}]}{T} = 0.
    \end{aligned}
\end{equation}

To solve the above problem, we define $\mathbf{{\Theta}}_t \stackrel{\triangle}{=} (\mathbf{J}_t,\mathbf{H}_t,Q_{c,t})$ as a concatenated vector of all queues, where $\mathbf{J}_t\stackrel{\triangle}{=} (J_{1,t},...,J_{K,t})$ and $\mathbf{H}_t \stackrel{\triangle}{=} (H_{1,t},...,H_{K,t})$.

The Lyapunov function, being a scalar measure of $\mathbf{{\Theta}}_t$, can be defined as follow:
\begin{equation}
    L(\mathbf{{\Theta}}_t) \stackrel{\triangle}{=} \frac{1}{2}\beta Q_{c,t}^2 + \frac{1}{2}\sum_{k \in \mathcal{K}}J_{k,t}^2 + \frac{1}{2}\sum_{k \in \mathcal{K}}H_{k,t}^2,
\end{equation}
where $\beta$ is a weight coefficient to treat the carbon-aware virtual queue differently from other queues such that we can adjust the emphasis placed on it. 

Then we define $\Delta (\mathbf{{\Theta}}_t) \stackrel{\triangle}{=} \mathbb{E}[L(\mathbf{{\Theta}}_{t+1}) - L(\mathbf{{\Theta}}_t) | \mathbf{{\Theta}}_t]$ as the conditional one-slot Lyapunov drift, measuring the expectation of the change in the Lyapunov function. Intuitively, a large value of $L(\mathbf{{\Theta}}_t)$ means that at least the backlog of one queue is large, which may stem from either substantial unfulfilled energy demand or higher carbon emission. Therefore, minimizing $\Delta (\mathbf{{\Theta}}_t)$ helps to push the queue state towards less congestion, which could maintain queue stability, but it may generate a low EV flexibility level. To overcome this issue, we introduce the drift-plus-penalty expression to balance the Lyapunov drift minimization and the EV flexibility maximization, i.e.,
\begin{equation}
 \min_{\check{p}_{k,t},\hat{p}_{k,t},\forall{k,t}}\; \Delta (\mathbf{{\Theta}}_t) + V(\check{p}_{s,t} -\hat{p}_{s,t})\Delta t,
\end{equation}
where $V$ is a positive parameter to adjust the impact of EV flexibility on the above objective function. For the sake of the time-coupling setting in $\Delta (\mathbf{{\Theta}}_t)$, we minimize the upper bound of the drift-plus-penalty term to implement the real-time algorithm. 
In the following, we will derive the upper bound of the drift-plus-penalty term.
Due to the fact that for any $q \geq 0, b \geq 0, a \geq 0$, we have:
\begin{equation} \label{eq:inequality_bound}
    (\max[q - b,0]+a)^2 \leq q^2 + a^2 + b^2 + 2q(a - b).
\end{equation}
Then based on \eqref{eq:inequality_bound}, we can have the following inequality for EV charging queue $J_{k,t}, k \in \mathcal{K}$:
\begin{equation}\label{J_k_bound}
    \begin{aligned}
        &\quad \;\frac{1}{2} \sum_{k \in \mathcal{K}} [J_{k,t+1}^2 - J_{k,t}^2]\\
        & =\frac{1}{2}\sum_{k \in \mathcal{K}}\left[(\max[J_{k,t} - \check{p}_{k,t},0]+a_{k,t})^2 - J_{k,t}^2\right]\\
        & \leq \frac{1}{2}\sum_{k \in \mathcal{K}} [a^2_{k,t} + \check{p}^2_{k,t}] + \sum_{k \in \mathcal{K}}J_{k,t}[a_{k,t} - \check{p}_{k,t}].
    \end{aligned}
\end{equation}

Using the fact that $H_{k+1,t} \leq \max[H_{k,t} + \frac{\lambda}{R_k} - \check{p}_{k,t},0]$, for delay-aware virtual queue we have:
\begin{equation}\label{H_k_bound}
    \begin{aligned}
    &\quad \;\frac{1}{2} \sum_{k \in \mathcal{K}} [H_{k,t+1}^2 - H_{k,t}^2]\\
    & \leq \frac{1}{2}\sum_{k \in \mathcal{K}} [\frac{\lambda}{R_k} - \check{p}_{k,t}]^2 + \sum_{k \in \mathcal{K}} H_{k,t}[\frac{\lambda}{R_k} - \check{p}_{k,t}].
    \end{aligned}
\end{equation}

Similarly, for the carbon-aware virtual queue, we have
\begin{equation}\label{Q_c_bound}
    \begin{aligned}
    &\quad \;\frac{1}{2} [Q_{c,t+1}^2 - Q_{c,t}^2]\\
    & \leq \frac{1}{2} [ w_{g,t}\hat{p}_{s,t} - r]^2 +  Q_{c,t}[ w_{g,t}\hat{p}_{s,t} - r].
    \end{aligned}
\end{equation}

Substituting \eqref{J_k_bound}, \eqref{H_k_bound} and \eqref{Q_c_bound} into the drift-plus-penalty expression yields:
\begin{equation}\label{all_bound}
    \begin{aligned}
         &\Delta (\mathbf{{\Theta}}_t) +  V (\check{p}_{s,t} -\hat{p}_{s,t})\Delta t  \leq  V (\check{p}_{s,t} -\hat{p}_{s,t})\Delta t\\
         & + \frac{1}{2}\sum_{k \in \mathcal{K}} [a^2_{k,t} + \check{p}^2_{k,t}] + \frac{1}{2}\sum_{k \in \mathcal{K}} [\frac{\lambda}{R_k} - \check{p}_{k,t}]^2 \\
         &+ \frac{1}{2} \beta (w_{g,t}\hat{p}_{s,t} - r)^2 + \sum_{k \in \mathcal{K}}J_{k,t} [a_{k,t} - \check{p}_{k,t}]\\
        &+ \sum_{k \in \mathcal{K}}H_{k,t} [\frac{\lambda}{R_k}-\check{p}_{k,t}] + \beta Q_{c,t} [w_{g,t}\hat{p}_{s,t} - r ].
    \end{aligned}
\end{equation}

\subsection{Real-Time Optimization Algorithm}
Minimizing the right-hand-side of \eqref{all_bound} leads to the following real-time optimization algorithm:

\begin{equation}\label{alg}
    \begin{aligned}
 \mathbf{P4\colon} &\min_{\check{p}_{k,t},\hat{p}_{k,t},\forall{k,t}}\;  \frac{1}{2}\beta \left(w_{g,t} \sum_{k \in \mathcal{K}}\hat{p}_{k,t} - r\right)^2 + \sum_{k \in \mathcal{K}}\check{p}_{k,t}^2\\
 & + \sum_{k \in \mathcal{K}}(V \Delta t - J_{k,t} - H_{k,t} - \frac{\lambda}{R_k})\check{p}_{k,t}\\
 & + \sum_{k \in \mathcal{K}}(-V \Delta t + \beta Q_{c,t}w_{g,t})\hat{p}_{k,t},\\
& \qquad \textrm{s.t.} \quad (\text{16d}) - (\text{16f}).
    \end{aligned}
\end{equation}

At each time slot, the aggregator would observe the real-time grid carbon-intensity $w_{g,t}$ and the current queue state $\mathbf{{\Theta}}_t$ to determine the real-time EV flexibility region $[\check{p}_{s,t}^\ast,\hat{p}_{s,t}^\ast]$, which is derived by solving $\mathbf{P4}$. We have $ \check{p}_{s,t}^\ast = \sum_{k \in \mathcal{K}}\check{p}_{k,t}^\ast$ and $ \hat{p}_{s,t}^\ast = \sum_{k \in \mathcal{K}}\hat{p}_{k,t}^\ast$, where $\check{p}_{k,t}^\ast$ and $\hat{p}_{k,t}^\ast$ are the optimal solutions of $\mathbf{P4}$.
Note that the proposed method is prediction-free. Then the solved flexibility interval is reported to the system operator for flexibility service provision. 

Another thing we care about is the gap between the result derived by the online algorithm $\mathbf{P4}$ and that of $\mathbf{P2}$, as shown in the following theorem:
\begin{theorem}\label{theorem:1}
        Suppose $a_{k,t}$, $w_{g,t}$, $\check{p}_{k,t}$, and $\hat{p}_{k,t}$ are upper bounded by $a_{k,max}$, $w_{g,max}$, $\check{p}_{k,max}$ and $\hat{p}_{k,max}$, respectively. Let $v_1^\ast$ and  $v^\ast$ denote the time-average power flexibility level expectation of the optimal solution obtained by $\mathbf{P2}$ and $\mathbf{P4}$, respectively. Then we have:
    \begin{equation}
        0 \leq v_1^\ast - v^\ast \leq \frac{B}{V},
    \end{equation}
where $B$ is a constant, defined as follows:
    \begin{equation}
    \begin{aligned}
        B \stackrel{\triangle}{=}&\frac{1}{2}\sum_{k \in \mathcal{K}} [a^2_{k,max} + \check{p}_{k,max}^2] +  \frac{1}{2}\sum_{k \in \mathcal{K}} \max[(\frac{\lambda}{R_k})^2 ,\check{p}_{k,max}^2] \\
        + &\frac{1}{2}\beta\max[(w_{g,max}\hat{p}_{s,max})^2,r^2].
    \end{aligned}
\end{equation}
\end{theorem}

The proof of Theorem \ref{theorem:1} is in Appendix C. Theorem \ref{theorem:1} tells us that a larger $V$ results in a narrower optimality gap, which, however, could increase the size of queues. On the contrary, a smaller value of $V$ leads to a larger optimality gap but gives a more stable queue state.

It is worth noting that different from existing works that design the real-time strategy by minimizing a linear function ($B + V(\check{p}_{s,t}^w -\hat{p}_{s,t}^w)\Delta t + \sum_{k \in \mathcal{K}}J_{k,t} [a_{k,t} - \check{p}_{k,t}^w]
        + \sum_{k \in \mathcal{K}}H_{k,t} [\frac{\lambda}{R_k}-\check{p}_{k,t}^w] + \beta Q_{c,t} [w_{g,t}\hat{p}_{s,t}^w - r ]$), the real-time algorithm \eqref{alg} minimizes a drift-plus-penalty expression by considering the quadratic terms, which can improve the performance of the online algorithm since the right-hand-side of \eqref{all_bound} is more tightly bounded by the quadratic function.
        
\emph{Remark:} It would be possible to explore how to integrate the available predictive information into the proposed method to further improve the performance of our algorithm.
One possible approach could involve incorporating a Lyapunov drift spanning multiple time slots (e.g., 2 time slots). This would require determining future queue backlogs, which in turn would require predictive information such as EV charging load forecasting. However, this possibility would need a more refined design and rigorous theoretical analysis since incorporating  predictions into the online decision-making does not necessarily lead to a better performance \cite{6322266}. Such developments are nontrivial and will be addressed in future work.

\section{Disaggregation and Queue Update}
\label{sec:Disaggregation}
\begin{figure}[!htpb]
    \centering
    \includegraphics[width = 0.9\linewidth]{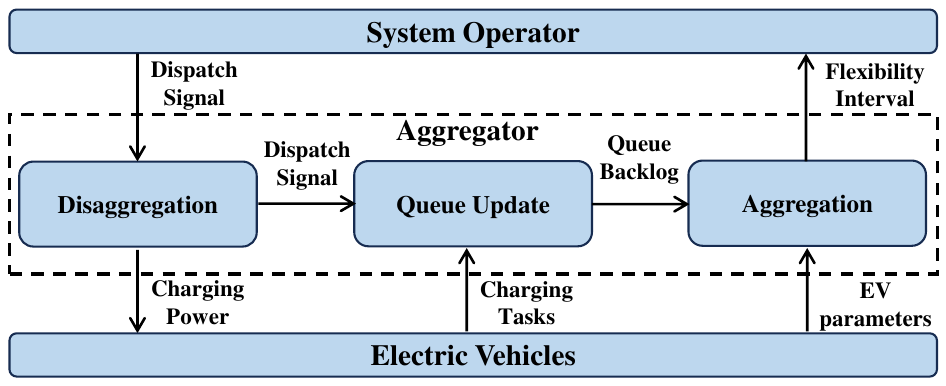} 
    \caption{Structure of the overall algorithm.}
    \label{fig:alg_structure}
\end{figure}

Based on the methods in Section \ref{sec:online}, the aggregator reports the flexibility interval $[\check{p}_{s,t}^\ast, \hat{p}_{s,t}^\ast]$ to the system operator in each time slot. 
Then the system operator can choose any power level within this flexibility interval and send a dispatch signal to the aggregator. The dispatch signal is disaggregated into the charging power of individual EVs by the aggregator to react to this order. In addiiton, when implementing the online algorithm \eqref{alg}, we conservatively use the upper bound $\hat{p}_{s,t}$ for the carbon-aware virtual queue update and the lower bound $\check{p}_{s,t}$ for the EV charging queue update since the dispatch signal is unknown. Therefore, after receiving the dispatch signal, it is natural for the aggregator to utilize this true value of the aggregate power to update the queue states before solving the problem for the next time slot because, intuitively, this would enhance EV flexibility. With this in mind, we perform the disaggregation process and update the queue using the dispatch signal before moving to the next time slot. Fig. \ref{fig:alg_structure} gives the overall structure of the proposed algorithm.

\subsection{Two-Stage Disaggregation}
Let $p_{s,t}^d$ be the dispatch signal at time slot $t$ such that $\check{p}_{s,t}^\ast \leq p_{s,t}^d \leq \hat{p}_{s,t}^\ast$. In this paper, we focus on the quantification of the EV flexibility on the aggregator side, and how the system operator determines the dispatch signal is ignored. 

We define a dispatch ratio $\gamma_t$ as follows:
\begin{equation}
    \gamma_t = \frac{p_{s,t}^d - \check{p}^\ast_{s,t}}{\hat{p}^\ast_{s,t} - \check{p}^\ast_{s,t}}.
\end{equation}

Then we decompose the dispatch signal $p_{s,t}^d$ into group-level dispatch signals ($p_{k,t}^d,\forall{k}$) with this dispatch ratio $\gamma_t$:

\begin{equation}
    p_{k,t}^d = (1 - \gamma_t)\check{p}_{k,t}^\ast + \gamma_t\hat{p}_{k,t}^\ast, \forall{k},
\end{equation}
where $p_{s,t}^d = \sum_{k \in \mathcal{K}}p_{k,t}^d$ is satisfied.

Next, the group-level dispatch signal $p_{k,t}^d$ will be allocated to the EVs in group $k$. Here we design a two-stage disaggregation process. The first stage aims to clear the charging tasks in EV charging queue backlog $J_{k,t}$ via a FIFO manner, where we use $p_{1,k,t}^d$ to denote the power level allocated in this stage. If there is still unallocated power remaining after the first stage, then the second stage comes in. At the second stage, the early arriving vehicle is given priority in power allocation until its power reaches its upper limit $\Bar{P}_{i,t}$. Then we continue to allocate the remaining power to the next early arriving EV until the remaining power is completely assigned. Let $p_{2,k,t}^d$ represent the allocated power in the second stage, and we have:
\begin{equation}
    p_{k,t}^d = p_{1,k,t}^d + p_{2,k,t}^d.
\end{equation}

\subsection{Queue Update Using Dispatch Signals}
After finishing the two-stage disaggregation process, we carry out the queue update using the obtained dispatch signals.

Since the first-stage disaggregation process is related to the EV charging queue backlog, we replace \eqref{eq:queue_J_k} for time slot $t$ with dynamics:
\begin{equation}
    J_{k,t+1} = \max[J_{k,t} - p_{1,k,t}^d,0]+a_{k,t}, \forall{k}.
\end{equation}
Similarly, for the delay-aware virtual queue, we replace \eqref{eq:queue_H_k} for time slot $t$ with dynamics as follows:
\begin{equation}
    H_{k,t+1} = \max[H_{k,t} + \frac{\lambda}{R_k} - p_{1,k,t}^d, 0], \forall{k}.
\end{equation}

For the carbon-aware virtual queue, we use the dispatch signal $p_{s,t}^d$ to update the queue backlog for time slot $t$, which could reflect the true carbon footprint, as shown in the following:
\begin{equation}
    Q_{c,t+1} = \max[Q_{c,t} + w_{g,t}p_{s,t}^d - r,0].
\end{equation}

With the updated values of $J_{k,t+1}$, $H_{k,t+1}$ and $Q_{c,t+1}$, then we can solve $\mathbf{P4}$ for next time slot $t+1$. In Algorithm 1, we give a complete description of the proposed method for carbon-aware aggregate EV power flexibility quantification.
\begin{algorithm}
\caption{Carbon-Aware Online Optimization for Real-Time Aggregate EV Power Flexibility Quantification}
\begin{algorithmic}
    \State \textbf{Initialization:}
    \State $t \gets 1$. The aggregator sets $J_{k,1} = 0$, $H_{k,1} = 0$, $Q_{c,1} = 0$ and sets the parameter $V$, $\lambda$ and $\beta$.
    \State \textbf{1. Aggregate Flexibility Quantification:} 
    \State - The aggregator collects the EV arriving charging tasks during time slot $t$ and stores them to EV charging queues according to their charging durations.
     \State - Observe real-time grid carbon intensity $w_{g,t}$ and solve $\mathbf{P4}$ to get the flexibility interval $ [\check{p}_{s,t}^\ast,\hat{p}_{s,t}^\ast]$ at time slot $t$.
    \State - Report $ [\check{p}_{s,t}^\ast,\hat{p}_{s,t}^\ast]$ to the system operator.
    \State \textbf{2. System-Wide Dispatch:}
    \State - The system operator receives the flexibility interval. 
    \State - Generate the dispatch signal $p_{s,t}^d$ within $ [\check{p}_{s,t}^\ast,\hat{p}_{s,t}^\ast]$.

    \State \textbf{3. Two-Stage Disaggregation:}
    \State The First-Stage Disaggregation:
    \State - The aggregator clears EV charging queue backlog $J_{k,t}$ in a FIFO manner to determine $p_{1,k,t}^d$.
     \State The Second-Stage Disaggregation:
    \State - Determine the remaining power $p_{2,k,t}^d$ and  clear it according to the order of EV arriving times. ($p_{2,k,t}^d = 0$ means that there is no need to enter the second-stage disaggregation)
    
    \State \textbf{4. Queue Update Using Dispatch Signals:}
    \State - The aggregator updates $J_{k,t+1}$ and $H_{k,t+1}$ using $p_{1,k,t}^d$.
    \State - Update $Q_{c,t+1}$ using $p_{s,t}^d$.
    \State - Let $t \gets t+1$ and go back to the Step \textbf{1}.
\end{algorithmic}
\end{algorithm}

\section{Numerical Simulation}
\label{sec:case}

This section evaluates the effectiveness of the proposed real-time algorithm. All the experiments are finished in the MATLAB 2020a platform with an Inter i5-11300H processor and 16 GB RAM. The Gurobi solver is used.

\subsection{Simulation Setup}
\begin{figure}[!htpb]
    \centering
    \includegraphics[width = 0.85\linewidth]{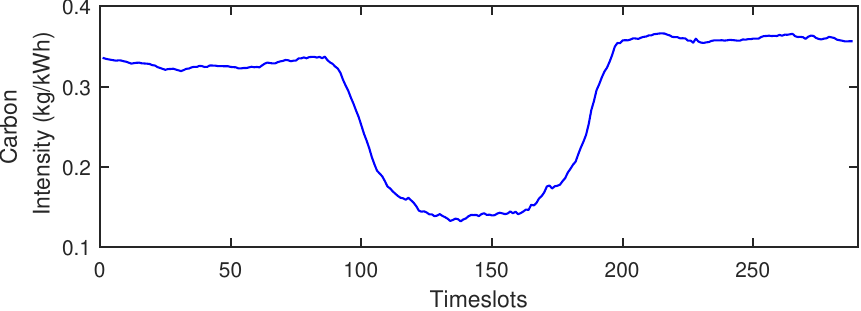} 
    \vspace{-1em}
    \caption{Grid carbon intensity profile.}
    \label{fig:carbon_data}
\end{figure}

The simulation horizon is 24 hours with 288 time slots. Each time slot has a time interval of 5 minutes. We use the real-world grid carbon intensity data from \cite{9960988}, whose profile is shown in Fig. \ref{fig:carbon_data}. 100 EVs are considered and we use normal distributions \cite{7790841} to describe the uncertainties of the EV charging parameters, which are presented in TABLE \ref{tab:EV para}. Based on the arrival and departure times of EVs, we can calculate the charging durations of all EVs, ranging from 4 hours to 12 hours. Therefore, we set $K$ = 9 groups and each group differs in one hour. The required SoC (the ratio between EV battery energy level and battery capacity) and the maximum SoC are set to 0.7 and 0.9, respectively. The charging efficiency is 0.95. We also take into account the diversity of the EVs with 3 types of EV battery capacities (60kWh, 40kWh, and 24kWh) and the corresponding maximum charging power (10kW, 6.6kW, and 3.3kW). We set $V$ to be 6000, $\beta$ to be 10 and $\lambda$ to be 100. The controlled carbon emission rate is chosen as $r$ = 30 kg/h.

\begin{table}[t]
\small
    \centering
    \caption{Parameters of the Stochastic EV charging} 
    \vspace{-0.5em}
    \begin{tabular}{@{}ccccccccc@{}}
        \toprule
        Parameter & Mean & Standard deviation\\
        \midrule
        Arrival time (h)   & 9 & 1.2  \\
        Departure time (h) & 18 & 1.2  \\
        Initial SoC & 0.4 & 0.1 \\
        \bottomrule
        \label{tab:EV para}
   \end{tabular}
\end{table}

\subsection{Benchmarks}
We use the following benchmarks to validate the performance of the proposed real-time method.

\begin{itemize}
    \item B1: \textit{Simple Algorithm}. In this method, each EV is charged at its maximum charging power $P_{i,max}$ upon arrival until its SoC reaches the required value of 0.7. Then the EV begins to show its flexibility with the lower bound of zero and upper bound of $P_{i,max}$. When the EV SoC reaches 0.9, the EV cannot provide any flexibility. Note that in each time slot, the carbon emission rate generated by the upper and lower bounds of the aggregate power of all EVs in service is strictly limited to no more than $r$ = 30 kg/h. (i.e., $w_{g,t}\check{p}_{s,t} \leq 30$, $w_{g,t}\hat{p}_{s,t} \leq 30$)
    
    \item B2: \textit{Greedy Algorithm}. This method can utilize the dispatch signal from the system operator. It first greedily chooses $a_t$ as the lower bound of the flexibility to avoid charging delay and $\sum_{i \in \mathcal{I}}\bar{P}_{i,t}$ as the upper bound to maximize the flexibility level. Similar to B1, the constraint of carbon emission rate is considered to limit the flexibility region but with a different strategy. To be specific, B2 greedily maintains the time-average carbon emission rate to be less than $r$ = 30kg/h, and thus the dispatch signal $p_{s,t}^d$ is used to update the maximum allowed carbon emission rate for each time slot.

    \item OPI: \textit{Offline Algorithm with Perfect Information}. This is the offline method by solving $\mathbf{P1}$ with the time-average carbon emission rate across all time slots less than $r$ = 30kg/h. Although it is impractical, this can provide a theoretical reference.

    \item B3: \textit{Real-Time Algorithm Based on Traditional Lyapunov Optimization}. B3 minimizes the linear expression of the drift-plus-penalty term rather than the quadratic one in the proposed method. Note that for a more rational comparison, B3 uses the same values for parameters $V$, $\beta$ and $\lambda$ as the proposed one.
\end{itemize}

\subsection{Performance Analysis}
\begin{figure}[!htpb]
    \centering
    \includegraphics[width = 0.85\linewidth]{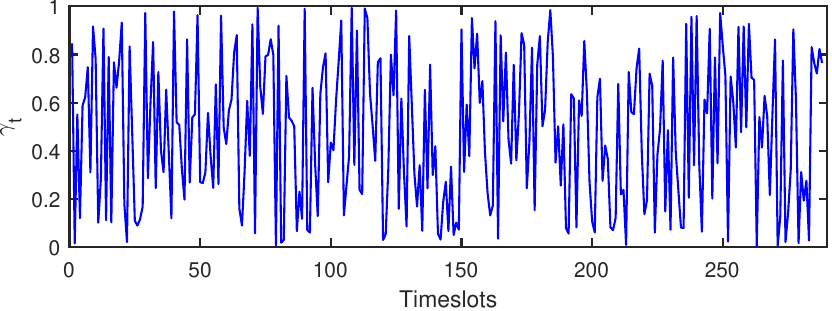} 
    \caption{Dispatch radio profile.}
    \label{fig:disp_ratio}
\end{figure}

\begin{figure}[!htpb]
    \centering
    \includegraphics[width = 0.95\linewidth]{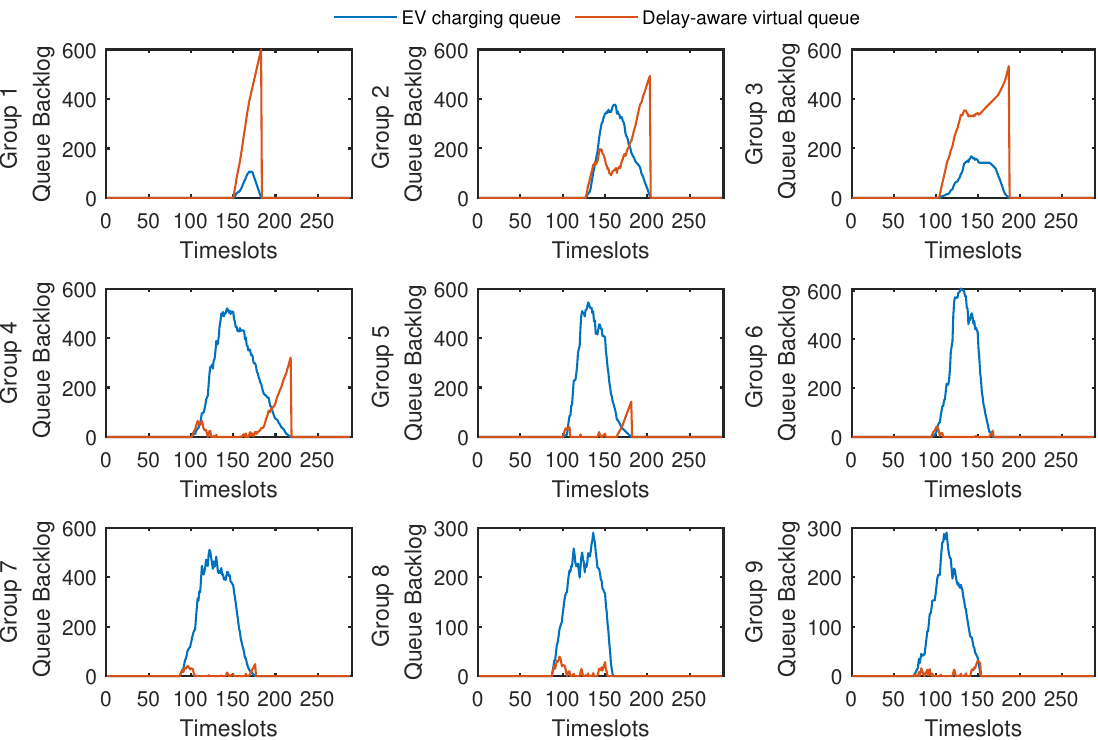} 
    \caption{Queue backlog profiles.}
    \label{fig:queue_backlog}
\end{figure}

Since the determination of the dispatch signal is out of the scope of this work, we randomly generate the dispatch ratio $\gamma_t \in [0,1]$ over the optimization horizon to conduct the case study, which is shown in Fig. \ref{fig:disp_ratio}. We first check the evolving queue backlogs of $J_{k,t}$ and $H_{k,t}$ in each group. In Fig. \ref{fig:queue_backlog}, we can see that the EV charging queue $J_{k,t}$ in each group first increases and then decreases, meaning the algorithm can utilize the temporal flexibility of EV charging behaviors by queue backlogs. We can observe the backlog of $J_{k,t}, \forall{k}$ converges to zero, which indicates all charging tasks of each group are fulfilled before the EVs depart. The delay awareness of the queue $H_{k,t},\forall{k}$ is reflected in the curves. One can observe that the groups (e.g., groups 1,2, and 3) with shorter charging durations tend to have a higher level of $H_{k,t}$ to push their charging tasks to be completed as soon as possible. However, the groups (e.g., groups 6, 7, and 8) with longer charging durations have a lower level of $H_{k,t}$, indicating these groups are almost not penalized by the $H_{k,t}$.

\begin{figure}[!htpb]
    \centering
    \includegraphics[width = 0.85\linewidth]{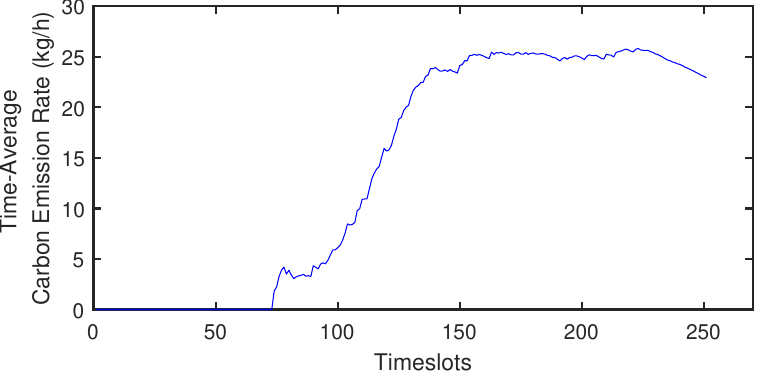} 
    \caption{Time-average carbon emission rate by the proposed method.}
    \label{fig:carbon_rate}
\end{figure}

Fig. \ref{fig:carbon_rate} shows the time-average carbon emission rate of the proposed online algorithm for each time slot. We can see that the time-average carbon emission rate converges to about 25 kg/h after a period of increase, and never exceeds the requirement ($r$ = 30 kg/h), showing the excellent control performance of the proposed algorithm on carbon emission.

\begin{figure}[!htpb]
    \centering
    \includegraphics[width = 0.95\linewidth]{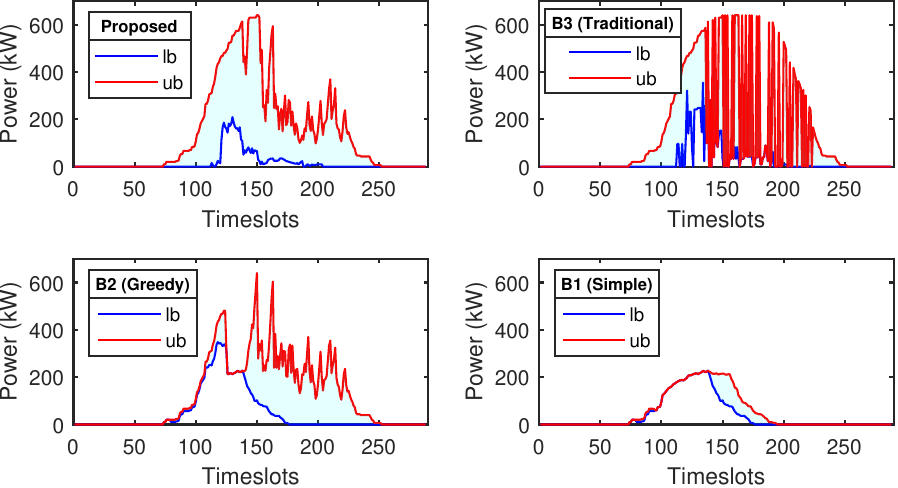} 
    \vspace{-1em}
    \caption{Flexibility regions by different methods.}
    \label{fig:flex_region}
\end{figure}

\begin{figure}[!htpb]
    \centering
    \includegraphics[width = 0.85\linewidth]{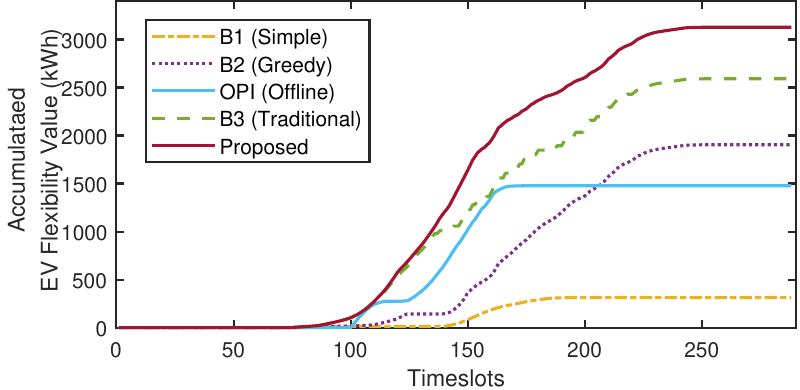} 
    \vspace{-1em}
    \caption{Time accumulated EV flexibility by different methods.}
    \label{fig:accu_flex}
\end{figure}

Fig. \ref{fig:flex_region} displays the aggregate EV flexibility regions of both the proposed method and benchmarks. Fig. \ref{fig:accu_flex} shows the accumulated EV flexibility value of all methods. TABLE \ref{tab:method_comp} compares four metrics of the methods involved, where the performance ratio is defined as the ratio of the total flexibility of a given method to the offline method. From Fig. \ref{fig:flex_region}, Fig. \ref{fig:accu_flex} and TABLE \ref{tab:method_comp}, we can see that the proposed algorithm achieves the maximum total flexibility and the largest performance ratio, which implies the advantages of our method over the traditional Lyapunov optimization (B3). Moreover, the proposed method shows a more evenly distributed flexibility region than B3 whose flexibility level oscillates up and down. B3 achieves a larger performance ratio than B2 since the greedy algorithm cannot effectively utilize the charging delay of EV tasks. Therefore, to greedily finish the charging tasks, B2 shows almost no significant flexibility before the time index of 150. However, the performance ratio of B2 is still greater than 1 because it accounts for the dispatch signal for update. The simple method B1 has the smallest total flexibility with a performance ratio of 0.21. Regarding the time-average carbon emission rate, all algorithms involved can keep it below the required value of 30. The proposed method and B3 perform best with the lowest carbon emission rate. Furthermore, the proposed method enables all required charging tasks to be fulfilled, while the simple method has a lot of unfilled energy, since it does not take into account the dispatch signal. B2 and B3  also have a small amount of unfulfilled energy.

\begin{table}[t]
\small
    \centering
    \caption{Results under different methods} 
    \vspace{-0.5em}
    \begin{tabular}{@{}ccccccccc@{}}
        \toprule
        Method & Proposed & B1 & B2 & B3\\
        \midrule
        Emission rate (kg/h)  & 22.94 & 12.15 & 26.29 & 22.93\\ Unfulfilled energy (kWh)  & 0 & 218.6 & 2.2 & 10.0\\ 
        Total flexibility (kWh) & 3127 & 312.3 & 1904 & 2590\\
        Performance ratio & 2.12 & 0.21 & 1.29 & 1.75\\
        \bottomrule
        \label{tab:method_comp}
   \end{tabular}
\end{table}

\begin{figure}[!htpb]
    \centering
    \includegraphics[width = 0.85\linewidth]{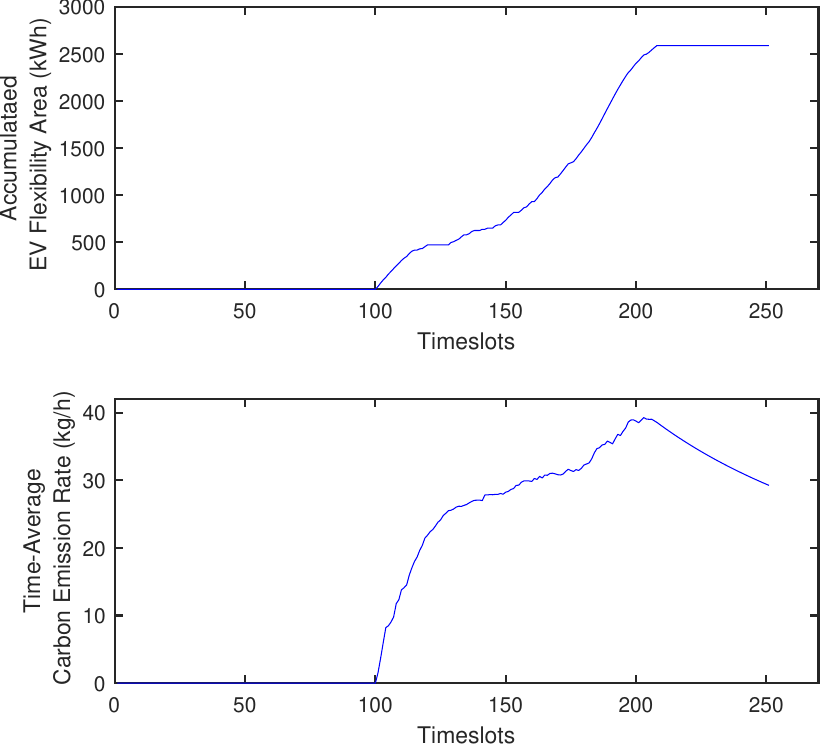} 
    \caption{Time accumulated EV flexibility and time-average carbon emission rate by the MPC method.}
    \label{fig:mpc}
\end{figure}

To further show the advantages of our method, we compare the proposed method with the MPC scheme. In the MPC setting, it is assumed that perfect knowledge of future EV arrival and grid carbon intensity is known. At each time slot, the aggregator solves the offline model \eqref{eq:off_obj} to maximize the total flexibility for future time slots, but only the flexibility interval of the first time slot is reported to the system operator. Upon receiving the dispatch signal from the system operator, the aggregator adopts the disaggregation method specified in Appendix A to ensure the feasibility of the problem. After the disaggregation process, the aggregator updates the states of all in-station EVs and the remaining carbon quota using the dispatch signal. At the next time slot, the aggregator repeats the same receding horizon optimization procedure. Fig \ref{fig:mpc} shows the results of the MPC method. The total EV flexibility of the MPC is 2589 kWh and the time-average carbon emission rate is 29.27 kg/h. Although MPC has perfect predictions of future uncertainties, it generates a lower total flexibility value and a higher carbon emission rate than our approach. It is worth noting that our approach demonstrates advantages over the MPC even without any prediction information. As a result, our Lyapunov optimization-based approach is more suitable for quantifying the real-time aggregate EV flexibility with carbon emission control, as our novel queue system accurately captures the structure of the problem with theoretical guarantees.

\subsection{Impact of Parameters}
\begin{figure}[!htpb]
    \centering
    \includegraphics[width = 0.95\linewidth]{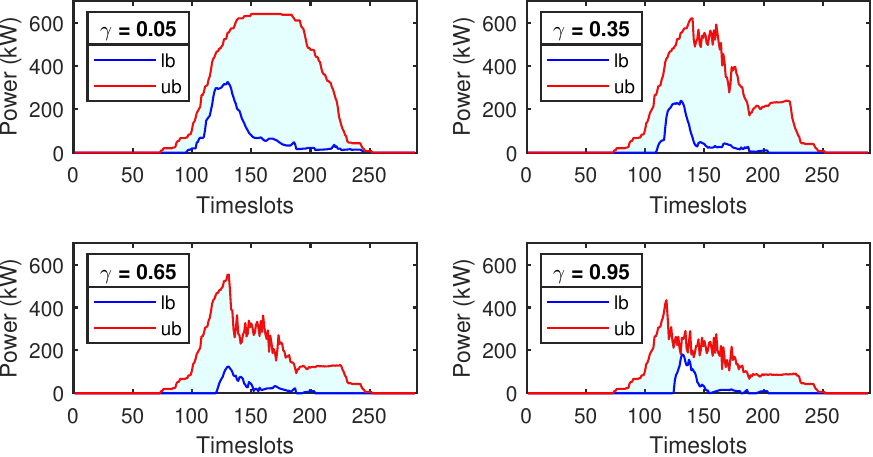} 
    \vspace{-1em}
    \caption{Flexibility regions by the proposed method with different $\gamma$.}
    \label{fig:gamma_region}
\end{figure}

We first discuss the effect of dispatch ratio $\gamma_t$ on the flexibility region. For ease of comparison, we use the fixed dispatch radio across all time slots (i.e., $\gamma_t = \gamma, \forall{t}$) instead of randomly generating it. Fig. \ref{fig:gamma_region}
presents the aggregate power flexibility regions under different $\gamma$ and TABLE \ref{tab:gamma_region}
shows the corresponding four metrics. We can see that a smaller value of $\gamma$ leads to a larger flexibility region, thus a greater performance ratio. As can be observed, when $\gamma$ increases, the power levels of both the lower and upper bounds of the flexibility region decrease to adapt to the change of $\gamma$. Specifically, a large value of $\gamma$ means that the system operators tend to dispatch much more charging power, which intuitively accelerates the accumulation of carbon emission, and in return, the charging tasks can be completed faster. The proposed algorithm adaptively reduces the power level of the upper bound to control the carbon emission rate, and simultaneously cuts down the power level of the lower bound due to the low intensity of the charging tasks. Moreover, in TABLE \ref{tab:gamma_region}, although the carbon emission rate grows as $\gamma$ increases, the rate is always controlled to be less than the required value of 30. Besides, we can see that the proposed method can almost guarantee the fulfillment of all charging tasks under any level of $\gamma$. Despite the 18.4 kWh of unfilled energy with $\gamma$ = 0.05, this amount of energy can be ignored with regard to the charging demand of 100 EVs.

\begin{table}[t]
\small
    \centering
    \caption{Results under the proposed method with different $\gamma$} 
    \vspace{-0.5em}
    \begin{tabular}{@{}ccccccccc@{}}
        \toprule
        $\gamma$ & 0.05 & 0.35 & 0.65 & 0.95\\
        \midrule
        Emission rate (kg/h)  & 14.86 & 21.82 & 24.38 & 25.88\\ Unfulfilled energy (kWh)  & 18.40 & 0 & 0 & 0\\ 
        Total flexibility (kWh) & 4363.0 & 3428.5 & 2496.4 & 1748.6\\
        Performance ratio & 2.95 & 2.32 & 1.69 & 1.18\\
        \bottomrule
        \label{tab:gamma_region}
   \end{tabular}
\end{table}

\begin{figure}[!htpb]
    \centering
    \includegraphics[width = 0.85\linewidth]{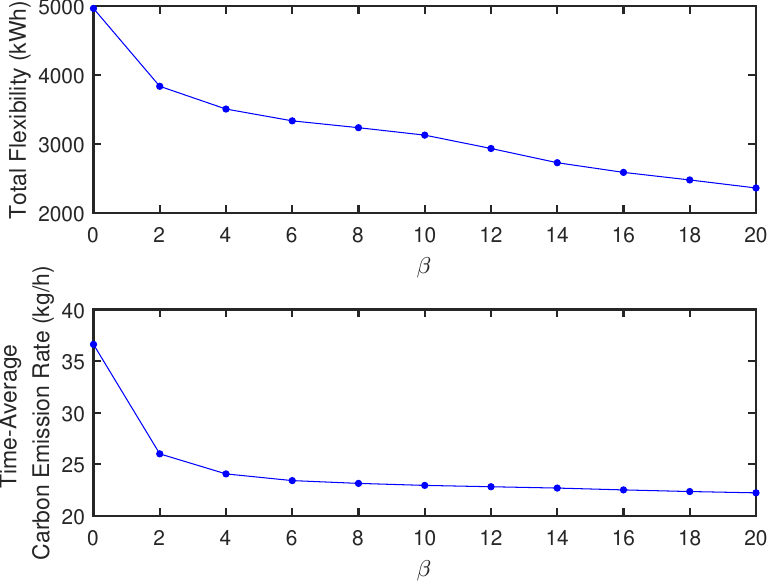} 
    \vspace{-1em}
    \caption{Total EV flexibility and time-average carbon emission rate by the proposed method with different $\beta$.}
    \label{fig:disc_beta}
\end{figure}

Then we investigate the impact of $\beta$ (the weight parameter of the carbon-aware virtual queue). We change $\beta$ from 0 to 20 and the results are shown in Fig. \ref{fig:disc_beta}. It can be observed that a higher value of $\beta$ results in both a lower total flexibility and a lower average carbon emission rate over time. This is because more carbon awareness is integrated into the algorithm by putting more emphasis on carbon-aware virtual queue stability, and the performance indicated by the total flexibility would decrease accordingly due to the reduction of carbon emissions.

\begin{figure}[!htpb]
    \centering
    \includegraphics[width = 0.85\linewidth]{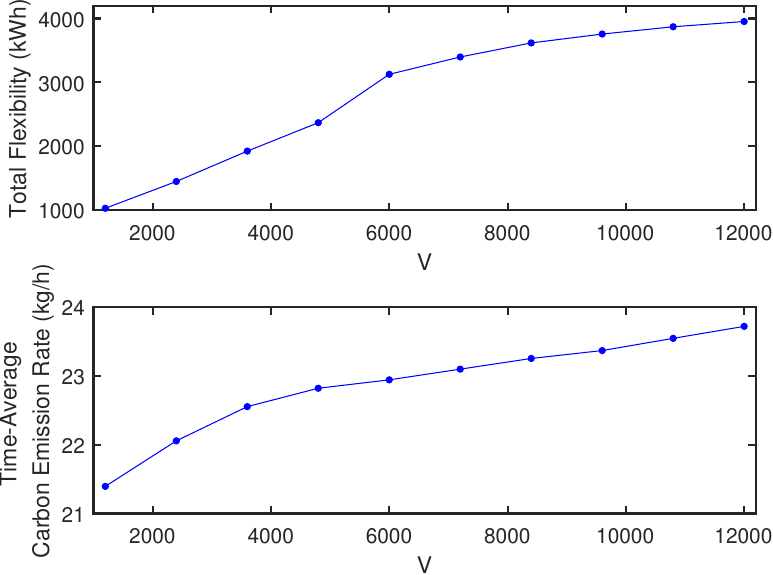} 
    \vspace{-1em}
    \caption{Total EV flexibility and time-average carbon emission rate by the proposed method with different $V$.}
    \label{fig:V_disc}
\end{figure}

We also analyze the results under different $V$, as shown in Fig. \ref{fig:V_disc}. We change the value of $V$ from 1200 to 12000. Since $V$ controls the trade-off between maximizing the EV flexibility level and stabilizing the queues, a large value of $V$ would contribute to a large time-average carbon emission rate, but incur greater total flexibility. This is consistent with the theoretical analysis in Theorem 1.

When the required energy of each EV is equal to the upper bound of its battery energy level (i.e., $e_i^{req} = e_{i, max}$), the offline model may fail to generate the EV flexibility region since the upper and lower trajectories of the EV flexibility region coincide. To confirm that our algorithm can still effectively quantify the real-time flexibility intervals under this situation, we set the parameter $e_{i}^{req}$ to the maximum SoC value of 0.9. TABLE \ref{tab:method_max_req} shows the corresponding results under different online benchmark algorithms. Our approach outperforms other algorithms, generating maximum total EV flexibility while keeping the lowest carbon emissions rate. We also investigate each algorithm's fulfillment ratio (i.e., the ratio of fulfilled charging requests to the total charging demand of all EVs) as the required energy level of EVs increases to 0.9. One can observe that our algorithm achieves the highest fulfillment ratio, which demonstrates its advantages.

\begin{table}[t]
\small
    \centering
    \caption{Results under different methods when the required energy of all EVs is the maximum battery energy limit}
    \vspace{-0.5em}
    \begin{tabular}{@{}ccccccccc@{}}
        \toprule
        Method & Proposed & B1 & B2 & B3\\
        \midrule
       Emission rate (kg/h)  & 25.79 & 13.92 & 26.34 & 26.91\\ 
        Total flexibility (kWh) & 2320.2 & 0 & 1145.0 & 1847.9\\
        Fulfillment ratio  & 0.974 & 0.603 & 0.928 & 0.968\\ 
        \bottomrule
        \label{tab:method_max_req}
   \end{tabular}
\end{table}

\begin{figure}[!htpb]
    \centering
    \includegraphics[width = 0.85\linewidth]{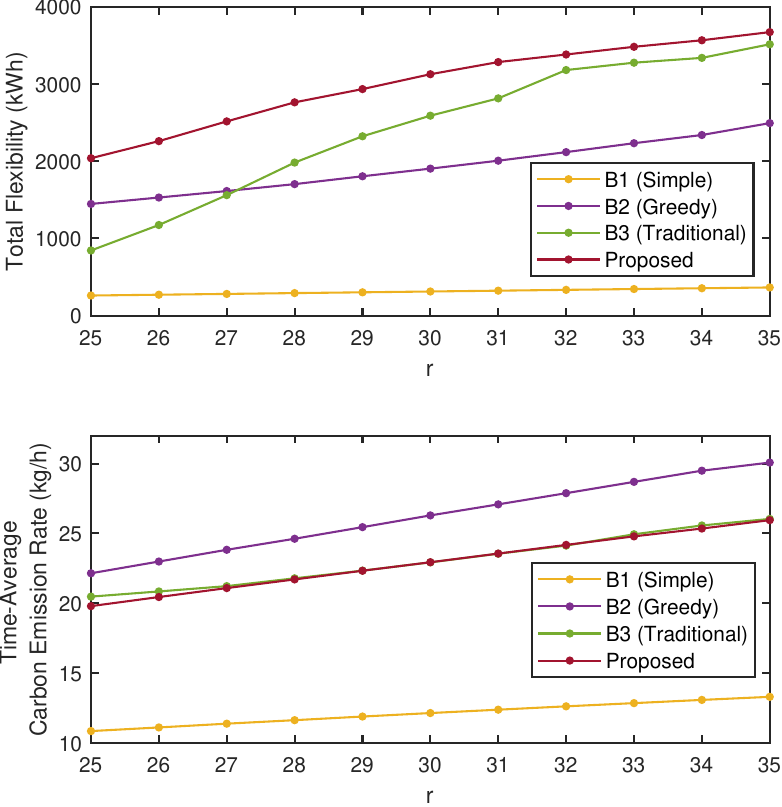} 
    \vspace{-1em}
    \caption{Total EV flexibility and time-average carbon emission rate by different methods with various $r$.}
    \label{fig:disc_carbon}
\end{figure}

Fig. \ref{fig:disc_carbon} illustrates the influence of carbon constraints (i.e., the value of $r$) on EV flexibility and carbon emission. As the value of $r$ increases from 25 to 35, the EV flexibility and carbon emissions of all algorithms increase due to the relaxation of the carbon emission constraint. Compared to the B3, the proposed method consistently exhibits greater flexibility with different $r$ while achieving nearly the same level of time-average carbon emission rate since it minimizes a quadratic form of the drift-plus-penalty term. In general, the B2 performs worse than B3, exhibiting lower flexibility and higher carbon emission. Interestingly, when $r$ is less than 27, one can observe that B2 achieves higher flexibility than B3. Combined with previous results in Fig. \ref{fig:flex_region} showing oscillation phenomenon in the performance of B3, this highlights the instability of B3 and further underscores the superiority of our proposed approach. B1 produces the least EV flexibility, resulting in the lowest carbon emission, but at the cost of a significant amount of unfulfilled energy.

\subsection{Scalability}

\begin{table}[t]
\small
    \centering
    \caption{Computational efficiency of offline optimization and proposed online optimziation}
    \vspace{-0.5em}
    \begin{tabular}{@{}ccccccccc@{}}
        \toprule
       Number of EVs & 100 & 200 & 500 & 1000\\
        \midrule
       Proposed & 0.0023s & 0.0025s & 0.0032s & 0.0040s\\ 
        Offline & 0.34s & 0.65s & 1.62s & 3.27s\\
        \bottomrule
        \label{tab:com_time}
   \end{tabular}
\end{table}

\begin{figure}[!htpb]
    \centering
    \includegraphics[width = 0.85\linewidth]{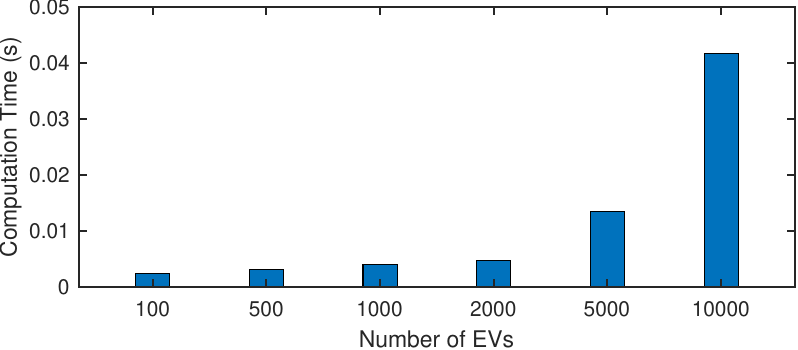} 
    \vspace{-1em}
    \caption{Average computation time of the proposed method for each charging decision with different numbers of EVs.}
    \label{fig:scalability}
\end{figure}

To show the practicability of the proposed method, we conduct a scalability analysis by increasing the number of EVs. Owing to the surge in charging demand driven by the large-scale growth in the number of EVs, we proportionally increased the cap of the time-average carbon emission rate $r$ to match the growth in charging demand. 
In TABLE VI, we first compare the computational efficiency of the offline model (i.e., OPI benchmark) and our proposed online model. When the number of EVs is 100, the average computation time of the online model for each charging decision is 0.0023s, while the offline model also achieves good computation efficiency with less than 0.5s for decision making. As the number of EVs gradually increases to 1,000, the computation times for both offline and online optimization models grow almost linearly. However, the online approach demonstrates significantly higher computational efficiency. Specifically, for every additional 100 EVs, the computation time of the online method increases only by around 0.0002s. In contrast, the offline model requires an additional 0.3s. This is because the number of variables in the offline model grows with the number of EVs, whereas our online model makes decisions at the queue level and the number of queues is independent of the number of EVs and bounded by $K$.
Then we further test the scalability of our proposed online method. Fig. \ref{fig:scalability} shows the average computation time per charging decision under different numbers of EVs. 
We observe that even when the number of EVs reaches 10,000, the algorithm can compute charging decisions within 0.05 seconds, which shows that our method is computationally efficient enough to handle a much larger number of EVs in a real-time setting and can be extended to the aggregate target of multiple charging stations in a microgrid since the network constraint is not time-coupled.


\section{Conclusion}
\label{sec:conclusion}
This paper develops a carbon-aware real-time algorithm for quantifying the aggregate EV power flexibility. The design of the real-time strategy is driven by the framework of Lyapunov optimization, which is prediction-free and can accommodate various uncertainties including EV arrivals and grid carbon intensity. The packaging technique is used to describe the EV charging tasks. Then, multiple virtual queues are used to control the charging delay and the carbon emission rate. The simulations validate the effectiveness of the proposed method and lead to the following findings:
\begin{enumerate}
    \item Compared with the existing real-time strategies, the proposed algorithm maximizes total EV flexibility while achieving a lower time-average carbon emission rate and zero unfulfilled energy.
    
   \item By minimizing a drift-plus-penalty expression with quadratic terms, the proposed method generates a larger and more uniform flexibility region than the traditional Lyapunov method.
    
    \item The dispatch signal affects the flexibility region. Specifically, increasing the dispatch ratio across all time slots can reduce the power levels of both upper and lower bounds.

\end{enumerate}

Future work will explore bi-directional (i.e., vehicle-to-grid) EV charging to enhance aggregate EV power flexibility.


\renewcommand\theequation{\Alph{section}.\arabic{equation}} 
\counterwithin*{equation}{section} 
\renewcommand\thefigure{\Alph{section}\arabic{figure}} 
\counterwithin*{figure}{section} 
\renewcommand\thetable{\Alph{section}\arabic{table}} 
\counterwithin*{table}{section} 

\begin{appendices}
\section{Proof of proposition 1}
 Let $p_{s,t}^{reg}$ be the EV aggregate regulation power at time $t$ such that $\check{p}^{\ast}_{s,t} \leq p_{s,t}^{reg} \leq \hat{p}^{\ast}_{s,t}$. We first determine an
auxiliary coefficient $\alpha_{t} \in [0,1]$ by:
\begin{equation}
    \alpha_{t} = \frac{\hat{p}^{\ast}_{s,t} -p_{s,t}^{reg}}{\hat{p}^{\ast}_{s,t} - \check{p}^{\ast}_{s,t}},
\end{equation}
which is followed by an observation that $p_{s,t}^{reg} = \alpha_{t}\check{p}^{\ast}_{d,t}+(1-\alpha_{t})\hat{p}^{\ast}_{d,t}$. Then, a feasible EV dispatch strategy for all time slots $t \in \mathcal{T}$ can be constructed by the following convex combinations:
\begin{subequations}
    \begin{align}
    p_{i,t} &= \alpha_{t}\check{p}^{\ast}_{i,t} + (1-\alpha_{t})\hat{p}^{\ast}_{i,t}, \\
    e_{i,t} &= \alpha_{t}\check{e}^{\ast}_{i,t} + (1-\alpha_{t})\hat{e}^{\ast}_{i,t} .
    \end{align}
\end{subequations}

We move on to prove that it is a feasible solution for  disaggregation as follows:
\begin{equation}
    \begin{aligned}
         p_{s,t}^{reg} &= \alpha_{t}\check{p}^{\ast}_{s,t}+(1-\alpha_{t})\hat{p}^{\ast}_{s,t} \\
         & = \alpha_{t}\sum_{i \in \mathcal{I}} \check{p}_{i,t}^{\ast} + (1-\alpha_{t})\sum_{i \in \mathcal{I}} \hat{p}_{i,t}^{\ast}\\
         & = \sum_{i \in \mathcal{I}} [\alpha_{t}\check{p}_{i,t}^{\ast} + (1-\alpha_{t})\hat{p}_{i,t}^{\ast}] = \sum_{i \in \mathcal{I}} p_{i,t},\\
    \end{aligned}
\end{equation}
which shows that constraint \eqref{eq:off_ub_1} holds for $p_{s,t}^{reg}$ and $p_{i,t}, \forall i$. Similarly, we can also prove that constraints \eqref{eq:off_ub_2}-\eqref{eq:off_ub_6} are satisfied. 

Since $p_{s,t}^{reg} \leq \hat{p}^{\ast}_{s,t}, \forall{t \in \mathcal{T}}$, we have
\begin{equation}
    \sum_{t \in \mathcal{T}}w_{g,t}p_{s,t}^{reg} \leq \sum_{t \in \mathcal{T}}w_{g,t}p^{\ast}_{s,t} \leq E_s,
\end{equation}
which shows that constraint \eqref{eq:off_carbon} is satisfied. 
Hence, we have found a carbon-aware feasible EV dispatch strategy for a given aggregate EV regulation power, which completes the proof.

\section{Proof of lemma 1}
Without loss of generality, we fix a time slot $t$. Then we demonstrate that for queue $k \in \mathcal{K}$, all EV charging task arrivals $a_{k,t}$ are finished at or before time slot $t + \delta_{k,max}$. If this is not true, we arrive at a contradiction. Observing from \eqref{eq:queue_J_k} that charging task arrivals are appended to the queue backlog $J_{k,t+1}$ and are first available for dispatch at time slot $t+1$. It has to be that $J_{k,\tau} > \check{p}_{k,\tau}$ for all $\tau \in \{t+1, ..., t+\delta_{k,max}\}$; otherwise the charging task backlog $a_{k,t}$ would be cleared at some time slot $\tau$ since the charging tasks are served in 
 FIFO order. Thus, by \eqref{eq:queue_H_k}, for all $\tau \in \{t+1, ..., t+\delta_{k,max}\}$, we have:
\begin{equation}
    H_{\tau,t+1} \geq H_{\tau,t} + \frac{\lambda}{R_k} - \check{p}_{k},
\end{equation}

Summing the above equation over $\tau \in \{t+1, ..., t+\delta_{k,max}\}$ gives:
\begin{equation}
    H_{k,t+\delta_{k,max} + 1} - H_{k,t+1} \geq \sum_{\tau = t+1}^{t+\delta_{k,max}}[-\check{p}_{k,\tau}] + \delta_{k,max}\frac{\lambda}{R_k},
\end{equation}

Since $H_{k,t+1} \geq 0$ and $H_{k,t+\delta_{k,max}+1} \leq H_{k,max}$, we derive the following inequality by rearranging the terms in (B.3):
\begin{equation}
    \delta_{k,max}\frac{\lambda}{R_k} \leq \sum_{\tau = t+1}^{t+\delta_{k,max}}[\check{p}_{k,\tau}] + H_{k,max},
\end{equation}

In particular, the sum of $\check{p}_{k,\tau}$ from slot $t+1$ to $t+\delta_{k,max}$ must be strictly less than $J_{k,t+1}$; otherwise, the charging task $a_{k,t}$, which is added at the end of the charging task backlog $J_{k,t+1}$, would be finished at a certain slot $\tau \in \{t+1, ..., t+\delta_{k,max}\}$ by the FIFO principle. Therefore:
\begin{equation}
    \sum_{\tau = t+1}^{t+\delta_{k,max}}[\check{p}_{k,\tau}] < J_{k,t+1} \leq J_{k,max}.
\end{equation}

Combining (B.4) and (B.3) leads to :
\begin{equation}
    \delta_{k,max}\frac{\lambda}{R_k} < J_{k,max} + H_{k,max},
\end{equation}
which indicates:
\begin{equation}
    \delta_{k,max} < \frac{(J_{k,max}+H_{k,max})R_k}{\lambda}.
\end{equation}

This contradicts the definition of $\delta_{k,max}$ in \eqref{def_delay}, proving that the worst-case charging task delay is bounded by $\delta_{k,max}$.

\section{Proof of theorem 1}
According to the theory of Lyapunov optimization\cite{neely2010stochastic}, for any given $\Tilde{\epsilon} > 0 $, there exists a so-called $w$-only policy $\check{p}_{k,t}^w, \hat{p}_{k,t}^w, \forall{t}$ that is feasible to \textbf{P2} with performance guarantee $\Tilde{\epsilon}$. Let $\Tilde{v}$ be the corresponding time-average flexibility energy level expectation. Furthermore, the $w$-only policy is independent of the virtual queue with $v_1^\ast \geq \Tilde{v} \geq v_1^\ast - \Tilde{\epsilon}$.

since $\check{x}_{k,t}^\ast, \hat{x}_{k,t}^\ast, \forall{t}$ is the optimal solution in \textbf{P4}, by \eqref{all_bound} we have:
\begin{align*}
         &\Delta (\mathbf{{\Theta}}_t^\ast) +  V (\check{p}_{s,t}^\ast -\hat{p}_{s,t}^\ast)\Delta t  \leq  V (\check{p}_{s,t}^w -\hat{p}_{s,t}^w)\Delta t\\
         & + \frac{1}{2}\sum_{k \in \mathcal{K}} [a^2_{k,t} + (\check{p}_{k,t}^w)^2] + \frac{1}{2}\sum_{k \in \mathcal{K}} [\frac{\lambda}{R_k} - \check{p}_{k,t}^w]^2 \\
         &+ \frac{1}{2} \beta (w_{g,t}\hat{p}_{s,t}^w - r)^2 + \sum_{k \in \mathcal{K}}J_{k,t} [a_{k,t} - \check{p}_{k,t}^w]\\
        &+ \sum_{k \in \mathcal{K}}H_{k,t} [\frac{\lambda}{R_k}-\check{p}_{k,t}^w] + \beta Q_{c,t} [w_{g,t}\hat{p}_{s,t}^w - r ] \\
        & \leq  V(\check{p}_{s,t}^w -\hat{p}_{s,t}^w)\Delta t  \tag{\stepcounter{equation}\theequation}\\
        &+ \frac{1}{2}\sum_{k \in \mathcal{K}} [a^2_{k,max} + \check{p}_{k,max}^2] +  \frac{1}{2}\sum_{k \in \mathcal{K}} \max[(\frac{\lambda}{R_k})^2 ,\check{p}_{k,max}^2] \\
        & + \frac{1}{2}\beta\max[(w_{g,max}\hat{p}_{s,max})^2,r^2] + \sum_{k \in \mathcal{K}}J_{k,t} [a_{k,t} - \check{p}_{k,t}^w] \\
        &  + \sum_{k \in \mathcal{K}}H_{k,t} [\frac{\lambda}{R_k}-\check{p}_{k,t}^w] + \beta Q_{c,t} [w_{g,t}\hat{p}_{s,t}^w - r ]\\
        & = B + V(\check{p}_{s,t}^w -\hat{p}_{s,t}^w)\Delta t + \sum_{k \in \mathcal{K}}J_{k,t} [a_{k,t} - \check{p}_{k,t}^w]\\
        &+ \sum_{k \in \mathcal{K}}H_{k,t} [\frac{\lambda}{R_k}-\check{p}_{k,t}^w] + \beta Q_{c,t} [w_{g,t}\hat{p}_{s,t}^w - r ]
\end{align*}

Utilizing the independence of the $w$-only policy on the virtual queue, we can get the following by taking expectations on both sides of (C.1):
\begin{equation}
    \begin{aligned}
        &\mathbb{E}[\Delta (\mathbf{{\Theta}}_t^\ast) +  V (\check{p}_{s,t}^\ast -\hat{p}_{s,t}^\ast)\Delta t ]\\
 & \leq B + V\mathbb{E}[(\check{p}_{s,t}^w -\hat{p}_{s,t}^w)\Delta t ]+\sum_{k \in \mathcal{K}}\mathbb{E}[J_{k,t} (a_{k,t} - \check{p}_{k,t}^w)]\\
        &+ \sum_{k \in \mathcal{K}}\mathbb{E}[H_{k,t} (\frac{\lambda}{R_k}-\check{p}_{k,t}^w)] + \beta \mathbb{E}[Q_{c,t} (w_{g,t}\hat{p}_{s,t}^w - r )]\\
        & = B + V\mathbb{E}[(\check{p}_{s,t}^w -\hat{p}_{s,t}^w)\Delta t ]+\sum_{k \in \mathcal{K}}\mathbb{E}[J_{k,t}] \mathbb{E}[a_{k,t} - \check{p}_{k,t}^w]\\
        &+ \sum_{k \in \mathcal{K}}\mathbb{E}[H_{k,t}] \mathbb{E}[\frac{\lambda}{R_k}-\check{p}_{k,t}^w] + \beta \mathbb{E}[Q_{c,t}] \mathbb{E}[w_{g,t}\hat{p}_{s,t}^w - r ]\\
        & \leq B + V\mathbb{E}[(\check{p}_{s,t}^w -\hat{p}_{s,t}^w)\Delta t ].
    \end{aligned}
\end{equation}

Summing the above equation over all time slots $t \in \{1,2,...,T\}$, dividing two sides by $VT$ and letting $T$ go to infinity , then we have:
\begin{equation}
    \begin{aligned}
       &\lim_{T \to \infty}\frac{1}{VT}\mathbb{E}[L(\mathbf{\Theta}_{T+1}^\ast)-L(\mathbf{\Theta}_{0}^\ast)] \\
       &+ \lim_{T \to \infty}\frac{1}{T}\sum_{t=1}^{T}\mathbb{E}[(\check{p}_{s,t}^\ast -\hat{p}_{s,t}^\ast)\Delta t]\\
        & \leq \frac{B}{V}+ \lim_{T \to \infty}\frac{1}{T}\sum_{t=1}^{T}\mathbb{E}[(\check{p}_{s,t}^w -\hat{p}_{s,t}^w)\Delta t].
    \end{aligned}
\end{equation}
Thus, we can obtain:
\begin{equation}
    \Tilde{v} \leq \frac{B}{V}+ v^\ast.
\end{equation}
Then we have:
\begin{equation}
    v_1^\ast - \Tilde{\epsilon} \leq \frac{B}{V} + v^\ast.
\end{equation}
Let $ \Tilde{\epsilon} \to 0$, we have:
\begin{equation}
     v_1^\ast - v^\ast\leq \frac{B}{V} .
\end{equation}
This completes the proof.
\end{appendices}

\bibliographystyle{IEEEtran}
\bibliography{references}

\end{document}